\documentclass[12pt,a4paper]{article}
\usepackage{setspace} 
\usepackage{amsfonts,amsmath,amssymb,amsthm,enumerate,CJK,array}
\usepackage{bm,dsfont,multicol}
\usepackage{graphicx,graphics,psfrag,booktabs}
\usepackage{natbib}
\usepackage[colorlinks=true]{hyperref} 
\usepackage{caption}
\usepackage{appendix}
\usepackage{titlesec}
\usepackage{rotating,threeparttable,booktabs}
\usepackage{longtable, lscape}
\usepackage{subfigure}

\DeclareMathOperator*{\argmin}{arg\,min}

\newcommand{\vep}{{\varepsilon}}

\newcommand{\CO}{{\mathcal{O}}}
\newcommand{\R}{\mathbb{R}}

\newcommand{\cov}{\mathop{\mbox{Cov}}}
\newcommand{\E}{\mbox{\sf E}}     
\renewcommand{\P}{\mathrm{P}}            

\newcommand{\IF}{\boldsymbol{1}}    
\def\Co{{\scriptstyle \mathcal{O}}} 
\def\defeq{\stackrel{\mathrm{def}}{=}}  

\def\CO{\mathcal{O}}
\def\Co{{\scriptstyle{\mathcal{O}}}}
\def\bx{{\mbox{\boldmath $x$}}}
\def\bX{{\mbox{\boldmath $X$}}}
\def\sbX{{\mbox{\scriptsize \boldmath $X$}}}
\def\bu{{\mbox{\boldmath $u$}}}
\def\sbu{{\mbox{\scriptsize \boldmath $u$}}}
\def\bv{{\mbox{\boldmath $v$}}}

\def\bt{{\mbox{\boldmath $t$}}}

\def\sbx{{\mbox{\scriptsize \boldmath $x$}}}

\def\boldalpha{{\mbox{\boldmath $\alpha$}}}
\def\sboldalpha{{\mbox{\scriptsize \boldmath $\alpha$}}}

\def\fx{f_{\mbox{\scriptsize \boldmath $X$}}}

\def\brh{\bar{h}}

\newtheorem{theorem}{THEOREM}[section]
\author{Shih-Kang Chao\thanks{Ladislaus von Bortkiewicz Chair of Statistics, C.A.S.E. - Center for applied Statistics and Economics, Humboldt-Universit\"{a}t zu Berlin, Unter den Linden 6, 10099
Berlin, Germany. email: shih-kang.chao@cms.hu-berlin.de; haerdle@wiwi.hu-berlin.de.} \and Katharina Proksch\thanks{Ruhr-Universit\"at Bochum, Fakult\"at f\"ur Mathematik, 44780 Bochum, Germany. email: katharina.proksch@rub.de; holger.dette@rub.de.}\and
 Holger Dette\footnotemark[3] \and Wolfgang
 H\"{a}rdle\footnotemark[2] \thanks{Lee Kong Chian School of Business, Singapore Management University, 50 Stamford Road, Singapore 178899, Singapore.} } 
\titleformat{\section}{\Large}{\textbf\thesection.}{.5em}{\textbf}
\titlespacing{\section}{0pt}{*3}{*2}
\titleformat{\subsection}{\large}{\textbf\thesubsection.}{.5em}{\textbf}
\titlespacing{\subsection} {0pt}{*3}{*2}
\titleformat{\subsubsection}{\normalfont}{\textbf\thesubsubsection.}{.5em}{\textbf}
\titlespacing{\subsubsection} {0pt}{*3}{*2}
\theoremstyle{definition} 

\newtheorem{corollary}[theorem]{Corollary}
\newtheorem{lemma}[theorem]{Lemma}

\newtheorem{remark}[theorem]{Remark}

\begin{document}
\title{Confidence Corridors for Multivariate Generalized Quantile Regression\thanks{Financial support from the Deutsche Forschungsgemeinschaft
(DFG) via SFB 649 "Economic Risk" (Teilprojekt B1), SFB 823
"Statistical modeling of nonlinear dynamic processes" (Teilprojekt
C1, C4) and Einstein Foundation Berlin via the Berlin Doctoral
Program in Economics and Management Science (BDPEMS) are gratefully
acknowledged.}} \maketitle
\begin{abstract}\singlespacing
We focus on the construction of confidence
corridors 
for multivariate nonparametric generalized quantile regression
functions. This construction is based on asymptotic results for the
maximal deviation between a suitable nonparametric estimator and the
true function of interest which follow after a series of
approximation steps including a Bahadur representation,  a new
strong approximation theorem and exponential tail inequalities for
Gaussian random fields.

As a byproduct we also obtain confidence corridors for the
regression function in the classical mean regression. In order to
deal with the problem of slowly decreasing error in coverage
probability of the asymptotic confidence corridors, which results in
meager coverage for small sample sizes, a simple bootstrap procedure
is designed based on the leading term of the Bahadur representation.
The finite sample properties of both procedures are investigated by
means of a simulation study and it is demonstrated that the
bootstrap procedure considerably outperforms the asymptotic bands in
terms of coverage accuracy. Finally, the bootstrap confidence
corridors are used to study the efficacy of the National Supported
Work Demonstration, which is a randomized employment enhancement
program launched in the 1970s. This article has supplementary
materials online.
\end{abstract}

\section{Introduction}\doublespacing
Mean regression analysis is a widely used tool in statistical
inference for curves. It focuses on the center of the conditional
distribution, 
given $d$-dimensional covariates with $d \geq 1$. In a variety of
applications though the interest is more in tail events, or even
tail event curves such as the conditional quantile function.
Applications with a specific demand in tail event curve analysis
include finance, climate analysis, labor economics and systemic risk
management.

Tail event curves have one thing in common: they describe the
likeliness of extreme events conditional on the covariate $\bX$. A traditional way
of defining such a tail event curve is by translating "likeliness"
with "probability" leading to conditional quantile curves. Extreme events
may alternatively be defined through conditional moment behaviour
leading to more general tail descriptions as studied by
\cite{NP:1987} and \cite{Jones:1994}. We employ this more general
definition
of generalized quantile regression (GQR),
which includes, for instance, expectile curves and study statistical inference of
GQR curves through confidence corridors.

In applications parametric forms are frequently used because of
practical numerical reasons. Efficient algorithms are available for
estimating the corresponding curves. However, the "monocular view"
of parametric inference has turned out to be too restrictive. This
observation prompts the necessity of checking the functional form of
GQR curves. Such a check may be based on testing different kinds of
variation between a hypothesized (parametric) model and a smooth
alternative GQR.
 Such an
approach though involves  either an explicit estimate of the bias or a
pre-smoothing of the "null model". In this paper we pursue the
Kolmogorov-Smirnov type of approach, that is,  employing the maximal deviation
between the null and the smooth GQR curve as a test statistic. Such a model check has
the advantage that it may be displayed graphically as a confidence
corridor (CC; also called
"simultaneous confidence band" or "uniform confidence band/region")
 but has been considered so far only for univariate covariates. The basic
technique for constructing  CC of this type is extreme value theory for the sup-norm
of an appropriately centered nonparametric estimate of the quantile curve.

Confidence corridors with one-dimensional predictor were developed under various settings. Classical one-dimensional results are confidence bands constructed for histogram estimators by \cite{smir:1950} or more general one-dimensional kernel density estimators by \cite{BR:1973}.
The results were extended to a univariate nonparametric mean regression setting by
\cite{Johnston82}, followed by \cite{WH:1989} who derived CCs for one-dimensional kernel $M$-estimators.
\cite{CK:2003} proposed  uniform confidence
bands and a bootstrap procedure for regression curves and their
derivatives.

In recent years, the growth of the literature body shows no sign of
decelerating. In the same spirit of \cite{WH:1989}, \cite{HS:2010}
and \cite{GH:2012} constructed  uniform confidence bands for
local constant quantile and expectile curves. \cite{FL:2013} proposed
an integrated approach for building simultaneous confidence band
that covers semiparametric models. \cite{GN:2010} investigated
adaptive density estimation based on linear wavelet and kernel
density estimators and \cite{LN:2011} extended the framework of
\cite{Bis:2007} to adaptive deconvolution density estimation.
Bootstrap procedures are proposed as a remedy for the poor coverage
performance of  asymptotic confidence corridors. For example, the
bootstrap for the  density estimator is proposed in \cite{H:1991}
and \cite{M:2012}, and for local constant quantile estimators in
\cite{SRH:2012}.

However, only recently progress has been achieved in the
construction of confidence bands for regression estimates with a
multivariate predictor. \cite{HH:12} derived an expansion for the
bootstrap bias and established a somewhat different way to construct
confidence bands without the use of extreme value theory. Their
bands are uniform with respect to a fixed but unspecified portion
(smaller than one) of points in a possibly multidimensional set in
contrast to the classical approach where uniformity is achieved on
the complete set considered. \cite{PBD:2014} proposed  multivariate
confidence bands for convolution type inverse regression models with
fixed design.

To the best of our knowledge, the classical Smirnov-Bickel-Rosenblatt type confidence corridors are not available for multivariate GQR or mean regression with random design.

In this work we go beyond the earlier studies in three aspects.
First, we extend the applicability of the CC to $d$-dimensional covariates with
 $d>1$. Second, we present a more general approach covering
not only quantile or mean curves but also GQR curves that are
defined via a minimum contrast principle. Third, we propose a
bootstrap procedure and we show numerically its improvement in the
coverage accuracy as compared to the asymptotic approach.

Our asymptotic results, which describe the maximal absolute
deviation of generalized quantile estimators, can not only be used
to derive a goodness-of-fit test in quantile  and expectile
regression, but they are also applicable in testing the quantile
treatment effect and stochastic dominance. We apply the new method
to test the quantile treatment effect of the National Supported Work
Demonstration program, which is a randomized employment enhancement
program launched in the 1970s. The data associated with the
participants of the program have been widely applied for treatment
effect research since the pioneering study of \cite{lalonde:1986}.
More recently, \cite{DE:2013} found that the program is beneficial
for individuals of over 21 years of age. In our study, we find that
the treatment tends to do better at raising the upper bounds of the
earnings growth than raising the lower bounds. In other words, the
program tends to increase the potential for high earnings growth but
does not reduce the risk of negative earnings growth. The finding is
particularly evident for those individuals who are older and spent
more years at school. We should note that the tests based on the
unconditional distribution cannot unveil the heterogeneity in
the earnings growth quantiles in treatment effects. 

The remaining part of this paper is organized as follows. In Section
\ref{Sec:Quantile} we present our model, describe the estimators and
state our asymptotic results. Section \ref{Sec:Bootstrap} is devoted
to the bootstrap and we discuss its theoretical and practical
aspects. The finite sample properties of both methods are
investigated by means of a simulation study in Section
\ref{Sec:simu}, where we also compare the numerical performance of our method with the method  proposed in \cite{HH:12} via simulations.
 The application of our new method is illustrated
by a real data example in Section \ref{Sec:apply}. The assumptions for
our asymptotic theory are listed and discussed after the references.
All detailed proofs are available in the supplement material.


\section{Asymptotic confidence corridors}\label{Sec:Quantile}
In Section \ref{pre} we present the prerequisites such as the
precise definition of the model and a suitable estimate. The result
on constructing confidence corridors (CCs) based on the distribution
of the maximal absolute deviation are given in Section
\ref{subsec:thms}. In Section \ref{subSec:scaling} we describe how
to estimate the scaling factors, which appear in the limit theorems,
using residual based estimators. Section \ref{subSec:bootstrap}
introduce a new bootstrap method for constructing CCs, while Section
\ref{subSec:bootQR} is devoted to specific issues related to
bootstrap CCs for quantile regression. Assumptions are listed and
discussed after the references.

\subsection{Prerequisites}\label{pre}
Let $(\bX_1,Y_1), ..., (\bX_n,Y_n)$ be a sequence of independent
identically distributed random vectors in $\Bbb R^{d+1}$ and
consider the nonparametric regression model
\begin{align}\label{model}
Y_i = \theta_0(\bX_i)+\varepsilon_i, \quad i=1,...,n,
\end{align}
where $\theta_0$ is an aspect of $Y$ conditional on $\bX$, such as the
$\tau$-quantile, the $\tau$-expectile or the mean regression curve, and the model errors $\varepsilon_1,\ldots,\varepsilon_n$ are i.i.d. with $\tau$-quantile,  $\tau$-expectile or  mean equal to 0, respectively, depending on which $\theta_0$ is in the model.
The function $\theta(\bx)$ can be estimated by:
\begin{align}
\hat \theta(\bx) = \argmin_{\theta \in \R}\frac{1}{n} \sum_{i=1}^n K_h(\bx-\bX_i) \rho(Y_i-\theta),
\end{align}
where $K_h(\bu)=h^{-d}K\left(\bu/h\right)$ for some kernel function
$K:\Bbb R^d \rightarrow \Bbb R$, and a loss-function $\rho_\tau:
\Bbb R \rightarrow \Bbb R$. 
In this paper we are concerned with the construction of uniform
confidence corridors for quantile as well as expectile regression
curves when the predictor is multivariate, that is, we focus on the
loss functions
$$\rho_\tau(u)=\big|\IF(u<0)-\tau\big||u|^k,$$
for $k=1$ and 2 associated with quantile and expectile regression.
We derive the asymptotic distribution of the properly scaled maximal
deviation $\sup_{\sbx\in\mathcal D}|\hat\theta_n(\bx)-\theta(\bx)|$
for both cases, where $\mathcal D \subset \R^d$ is a compact subset.
We use strong approximations of the empirical process, concentration
inequalities for general Gaussian random fields and results from
extreme value theory. To be precise, we show that
\begin{align}
&\P\left[(2 \delta \log n)^{1/2}\Big\{\sup_{\sbx \in \mathcal D}\big|r_n(\bx)\big[\hat \theta_n(\bx)-\theta_0(\bx)\big]\big|/\|K\|_2-d_n \Big\}<a \right] \rightarrow \exp\big\{-2 \exp(-a)\big\}, \label{max_dev}
\end{align}
as $n \rightarrow \infty$, where $r_n(\bx)$ is a scaling factor which
depends on $\bx$, $n$ and the loss function under consideration.
\subsection{Asymptotic results}\label{subsec:thms}
In this section we present our main theoretical results on the
distribution of the uniform maximal deviation of the quantile and
expectile estimator. The proofs of the theorems at their full
lengths are deferred to the appendix. Here we only give a brief
sketch of proof of Theorem \ref{TheoremQuantile} which is the limit
theorem for the case of quantile regression.

\begin{theorem}\label{TheoremQuantile}
Let $\hat \theta_n(\bx)$ and
$\theta_0(\bx)$ be the local constant quantile estimator and the
true quantile function, respectively and
suppose that assumptions (A1)-(A6) in Section \ref{Sec:Ass} hold. Let further $\mbox{vol}(\mathcal D)=1$ and
\begin{align*}
d_n = (2d \kappa \log n)^{1/2}+\big\{2d\kappa (\log n)\big\}^{-1/2} \left[\frac{1}{2}(d-1)\log \log n^\kappa + \log \big\{(2 \pi)^{-1/2} H_2 (2d)^{(d-1)/2}\big\} \right],
\end{align*}
where $ H_2 = \big(2\pi \|K\|_2^2\big)^{-d/2} \det(\Sigma)^{1/2},
\mbox{ }\Sigma=\bigl(\Sigma_{ij}\bigr)_{1 \leq i,j \leq d} =
\Bigl(\int \frac{\partial K(\sbu)}{\partial u_i } \frac{\partial
K(\sbu)}{\partial u_j } d\sbu\Bigr)_{1 \leq i,j \leq d}, $
\begin{align*}
r_n(\bx) = \sqrt{\frac{nh^d f_\sbX(\bx)}{\tau(1-\tau)}} f_{Y|\sbX}\big\{\theta_0(\bx)|\bx\big\},
\end{align*}
Then the limit theorem \eqref{max_dev} holds.
\end{theorem}

\noindent{\bf Sketch of proof.}
A major technical difficulty is imposed by the fact that the
loss-function $\rho_{\tau}$ is not smooth which means that standard
arguments such as those based on Taylor's theorem do not apply. As a consequence the use of a different, extended methodology
becomes necessary. In this context \cite{KLXbaha:2010} derived a
uniform Bahadur representation for an $M$-regression function in a
multivariate setting (see appendix). 
It holds uniformly for $x \in \mathcal D$, where $\mathcal D$ is a
compact subset of $\Bbb R^d$:
\begin{align}
\hat \theta_n(\bx)-\theta_0(\bx) = \frac{1}{n S_{n,0,0}(\bx)} \sum_{i=1}^n K_h(\bx-\bX_i)\psi_{\tau}\big\{Y_i-\theta_0(\bx)\big\} + \CO\left\{\Big(\frac{\log n}{nh^d}\Big)^{\frac{3}{4}}\right\}, \quad a.s. \label{baha1}
\end{align}
Here $S_{n,0,0}(\bx) = \int K(\bu) g(\bx+h \bu) f_\sbX(\bx+h\bu)
d\bu$, $\psi_\tau(u)=\IF(u<0)-\tau$ is the piecewise derivative of the loss-function $\rho_{\tau}$ and
\begin{align*}
g(\bx) &= \left.\frac{\partial}{\partial t} \E[\psi_{\tau}(Y-t)|\bX=\bx]\right|_{t=\theta_0(\sbx)}. 
\end{align*}

Notice that the error term of the Bahadur expansion does not depend
on the design $\bX$ and it converges to 0 with rate $\big(\log
n/nh^d\big)^{\frac{3}{4}}$ which is much faster than the convergence
rate $(nh^d)^{-\frac{1}{2}}$ of the stochastic term.

Rearranging \eqref{baha1}, we obtain
\begin{align}\label{emp1tilde}
S_{n,0,0}(\bx) \{\hat \theta_n(\bx)-\theta_0(\bx)\} = \frac{1}{n} \sum_{i=1}^n K_h(\bx-\bX_i)\psi_{\tau}\big\{Y_i-\theta_0(\bx)\big\}+\CO\left\{\Big(\frac{\log n}{nh^d}\Big)^{\frac{3}{4}}\right\}.
\end{align}
Now we express the leading term on the right hand side of \eqref{emp1tilde} by means of the centered empirical
process
\begin{align}
Z_n(y,\bu) = n^{1/2} \{F_n(y,\bu)-F(y,\bu)\},
\end{align}
where $F_n(y,\bx) = n^{-1} \sum_{i=1}^n \IF(Y_i \leq y, X_{i1} \leq
x_1, ..., X_{id} \leq x_d)$. This yields, by Fubini's theorem,
\begin{align}
S_{n,0,0}(\bx) \{\hat \theta_n(\bx)-\theta_0(\bx)\}-b(\bx) = n^{-1/2}\int \int K_h(\bx-\bu)\psi_{\tau}\big\{y-\theta_0(\bx)\big\}dZ_n(y,\bu)+\CO\left\{\Big(\frac{\log n}{nh^d}\Big)^{\frac{3}{4}}\right\}, \label{baha1.5}
\end{align}
where
$$
b(\bx) = -\E_\sbx\left[\frac{1}{n} \sum_{i=1}^n K_h(\bx-\bX_i)\psi\big\{Y_i-\theta_0(\bx)\big\}\right]
$$
denotes the bias which is of order $\CO(h^s)$ by Assumption (A3) in the Appendix. The variance of
the first term of the right hand side of \eqref{baha1.5} can be
estimated via a change of variables and Assumption (A5), which gives
\begin{align*}
&(nh^d)^{-2} n \E\big[K^2\big\{(\bx-\bX_i)/h\big\}\psi^2\big\{Y_i-\theta_0(\bx)\big\}\big]\\
&=(nh^d)^{-2} nh^d \int \int K^2(\bv) \psi^2\big\{y-\theta_0(\bx)\big\}f_{Y|\sbX}(y|\bx-h\bv) f_\sbX(\bx-h\bv)dyd\bv \\
&=(nh^d)^{-1} \int \int K^2(\bv) \psi^2\big\{y-\theta_0(\bx)\big\} f_{Y|\sbX}(y|\bx) f_\sbX(\bx)dyd\bv + \CO\big((nh^{d-1})^{-1}\big) \\
&=(nh^d)^{-1} f_\sbX(\bx) \sigma^2(\bx) \|K\|_2^2 + \CO\big\{(nh^{d})^{-1}h\big\},
\end{align*}
where $\sigma^2(\bx) =
\E[\psi^2\big\{Y-\theta_0(\bx)\big\}|\bX=\bx]$. The standardized
version of \eqref{emp1tilde} can therefore be approximated by
\begin{align}
&\frac{\sqrt{nh^d}}{\sqrt{f_\sbX(\bx) }\sigma(\bx)\|K\|_2}S_{n,0,0}(\bx) \{\hat \theta_n(\bx)-\theta_0(\bx)\}\notag \\
&= \frac{1}{\sqrt{h^d f_\sbX(\bx) }\sigma(\bx)\|K\|_2} \int \int K\left(\frac{\bx-\bu}{h}\right)\psi\big\{Y_i-\theta_0(\bx)\big\}dZ_n(y,\bu) + \CO\big(\sqrt{nh^d} h^s\big)+\CO\left\{\Big(\frac{\log n}{nh^d}\Big)^{\frac{3}{4}}\right\}.
\end{align}
The dominating term is defined by
\begin{align}\label{Y}
Y_n(\bx) \defeq \frac{1}{\sqrt{h^d f_\sbX(\bx)} \sigma(\bx)}\int \int K\left(\frac{\bx-\bu}{h}\right)\psi\big\{y-\theta_0(\bx)\big\} dZ_n(y,\bu).
\end{align}
Involving strong Gaussian approximation and
Bernstein-type concentration inequalities, this process can be
approximated by a stationary Gaussian field:
  \begin{align}\label{Y5}
Y_{5,n}(\bx)  = \frac{1}{\sqrt{h^{d}}} \int K\left(\frac{\bx-\bu}{h}\right) dW\big(\bu\big),
\end{align}
where $W$ denotes a Brownian sheet. The
 supremum  of this process is asymptotically Gumbel distributed, which follows, e.g., by Theorem 2 of \cite{Rosen:1976}.
Since the kernel is symmetric and
of order $s$, we can estimate the term
$$
S_{n,0,0} =  f_{Y|\sbX}(\theta_0(\bx)|\bx) f_\sbX(\bx) + \CO(h^s)
$$
if (A5) holds. On the other hand, $\sigma^2(\bx)=\tau(1-\tau)$ in
quantile regression. Therefore, the statements of the theorem hold.

\hfill$\Box$

\begin{corollary}[CC for multivariate quantile regression]\label{CorQuantile}
Under the assumptions of Theorem \ref{TheoremQuantile}, an
approximate $(1-\alpha)\times 100\%$ confidence corridor is given by
\begin{align*}
\hat \theta_n(\bt) \pm (nh^d)^{-1/2} \big\{\tau(1-\tau)\|K\|_2/\hat f_\sbX(\bt)\big\}^{1/2} \hat f_{\varepsilon|\sbX}\big\{0|\bt \big\}^{-1} \Big\{d_n+c(\alpha)(2 \kappa d \log n)^{-1/2} \Big\},
\end{align*}
where $\alpha
\in (0,1)$ and $c(\alpha) = \log 2- \log\big|\log(1-\alpha)\big|$ and $\hat
f_\sbX(\bt)$, $\hat f_{\varepsilon|\sbX}\big\{0|\bt\big\}$ are
consistent estimates for $\fx(\bt)$,
$f_{\varepsilon|\sbX}\big\{0|\bt\big\}$ with convergence rate in sup-norm faster
than $\Co_p\big((\log n)^{-1/2}\big)$.
\end{corollary}

\begin{remark}Note that under the conditions of Corollary \ref{CorQuantile} we find
	 \begin{align*}
		 \sup_{\sbx\in\mathcal{D}}\bigl|r_n(\bx)\bigl(\hat\theta_n(\bx)-\theta_0(\bx))\bigr)\bigr|=\CO_P\bigl(\sqrt{\log(n)}\bigr),
	\end{align*}
where 
	\begin{align*}
		 r_n(\bx)= \sqrt{\frac{nh^d f_\sbX(\bx)}{\tau(1-\tau)}} f_{Y|\sbX}\big\{\theta_0(\bx)|\bx\big\}.
  \end{align*}
For kernel estimators  $\hat f_{\varepsilon|\sbX}(0,\cdot)$	and $\hat f_{\sbX}(\cdot)$ converging in sup-norm with rate $\Co_P\bigl(\log(n)^{-1/2}\Bigr)$ to  $ f_{\varepsilon|\bX}(0,\cdot)$	and $f_{\bX}(\cdot)$, respectively, the quantity $\hat r_n(\bx)$, defined by
  \begin{align*}
	  \hat r_n(\bx)= \sqrt{\frac{nh^d \hat f_\sbX(\bx)}{\tau(1-\tau)}} \hat f_{\varepsilon|\sbX}(0,\bx),
	\end{align*}
inherits this rate. Furthermore, since we consider an additive error model, the conditional density $f_{Y|\sbX}\big\{\theta_0(\bx)|\bx\big\}$ can be replaced by $f_{\varepsilon|\sbX}(0,\bx)$ (see Section \ref{subSec:scaling} below for more details and the definition of suitable estimators).
This yields
  \begin{align*}
	   \sup_{\sbx\in\mathcal{D}}\bigl|\hat r_n(\bx)\bigl(\hat\theta_n(\bx)-\theta_0(\bx))\bigr)\bigr|=
		\Co_P(1)+\sup_{\sbx\in\mathcal{D}}\bigl|r_n(\bx)\bigl(\hat\theta_n(\bx)-\theta_0(\bx))\bigr)\bigr|.
	\end{align*}
	Hence, by Slutsky's Lemma,  the quantities $\sup_{\sbx\in\mathcal{D}}\bigl|\hat r_n(\bx)\bigl(\hat\theta_n(\bx)-\theta_0(\bx))\bigr)\bigr|$ and $\sup_{\sbx\in\mathcal{D}}\bigl| r_n(\bx)\bigl(\hat\theta_n(\bx)-\theta_0(\bx))\bigr)\bigr|$ have the same asymptotic distribution.
\end{remark}

The expectile confidence corridor can be constructed in an analogous manner  as  the quantile confidence corridor. The two cases differ in the
form  and hence the properties of the loss function. Therefore we find for expectile regression:
\begin{align*}
S_{n,0,0}(\bx) = -2 \big[F_{Y|\sbX}\big(\theta_0(\bx\big)|\bx)(2 \tau-1)-\tau\big]\fx(\bx) + \CO(h^s).
\end{align*}
Through similar approximation steps as the quantile regression, we
derive the following theorem.

\begin{theorem}\label{TheoremExpectile}
Let $\hat \theta_n(\bx)$  be the the local constant expectile estimator and $\theta_0(\bx)$ the
true expectile function.
If Assumptions (A1), (A3)-(A6) and (EA2) of Section \ref{Sec:Ass} hold with  a constant $b_1$ satisfying
$$
n^{-1/6} h^{-d/2-3d/(b_1-2)} = \CO(n^{-\nu}), \quad \nu >0.
$$
Then the limit theorem \eqref{max_dev} holds with a scaling factor $$r_n(\bx) = \sqrt{nh^d f_\sbX(\bx)}\sigma^{-1}(\bx)\left\{2 \big[\tau-F_{Y|\sbX}(\theta_0(\bx)|\bx)(2 \tau-1)\big]\right\}$$ and with the same constants $H_2$ and $d_n$ as defined in Theorem \ref{TheoremQuantile}, where $\sigma^2(\bx) = \E[\psi_{\tau}^2(Y-\theta_0(\bx))|\bX=\bx]$ and
$\psi_{\tau}(u)=2(\IF(u \leq 0)-\tau)|u|$ is the derivative of the expectile loss-function $\rho_{\tau}(u)=\big|\tau-\IF(u<0)\big| |u|^2$.
\end{theorem}
The proof of this result is deferred to the appendix. In the next
corollary,  the explicit form of the CCs for expectiles is given.
\begin{corollary}[CC for multivariate expectile regression] Under the same assumptions of
Theorem \ref{TheoremExpectile}, an approximate $(1-\alpha)\times
100\%$ confidence corridor is given by
\begin{align*}
\hat \theta_n(\bt) \pm (nh^d)^{-1/2} \big\{\hat \sigma^2(\bt)\|K\|_2/\hat f_\sbX(\bt)\big\}^{1/2} &\Big\{-2 \big[\hat F_{\varepsilon|\sbX}\big\{0|\bt \big\}(2 \tau-1)-\tau\big]\Big\}^{-1} \Big\{d_n+c(\alpha)(2 \kappa d \log n)^{-1/2} \Big\},
\end{align*}
where $\alpha \in (0,1)$ $c(\alpha) = \log 2- \log\big|\log(1-\alpha)\big|$ and $\hat
f_\sbX(\bt)$, $\hat \sigma^2(\bt)$ and $\hat
F_{\varepsilon|\sbX}(0|\bx)$ are consistent estimates for
$\fx(\bt)$, $\sigma^2(\bt)$ and $ F_{\varepsilon|\sbX}(0|\bx)$ with
convergence rate in sup-norm faster than $\Co_p\big((\log n)^{-1/2}\big)$.
\end{corollary}

A further immediate consequence of Theorem \ref{TheoremExpectile} is a similar limit theorem
in the context of local least squares estimation of the regression curve in classical mean regression.
\begin{corollary}[CC for multivariate mean regression]
Consider the loss function $\rho(u)=u^2$ corresponding to $\psi(u)=2u.$ Under the assumptions of Theorem \ref{TheoremExpectile}, with the same constants $H_2$ and $d_n$, \eqref{max_dev} holds for the local constant estimator $\hat\theta$ and the regression function $\theta(x)=E[Y\,|\,\bX=\bx]$ with scaling factor $r(\bx)=\sqrt{nh^d f_\sbX(\bx)}\sigma^{-1}(\bx)$ and  $\sigma^2(\bx)=$Var$[Y\,|\,\bX=\bx].$
\end{corollary}

\begin{remark}
We would like to stress that our purely non-parametric approach offers flexibility and reasonable results in moderate dimensions $d=2,\, d=3$, but it is not suitable for inference in high dimensional models due to the curse of dimensionality. The case of high dimensional regressors may be handled via a semi-parametric specification of the regression curve, such as, for instance, a partial linear model. Such a model was considered in \cite{SRH:2012} with a one-dimensional non-parametric component.
We think that our approach allows to adapt these ideas and, as an extension, to consider a non-parametric component which is multivariate. Hence, our approach then also offers higher flexibility in semi-parametric modeling. This semi-parametric approach is not pursued further in this paper but it clearly deserves future research.
\end{remark}

\subsection{Estimating the scaling factors}\label{subSec:scaling}

The performance of the confidence bands is greatly influenced by the
scaling factors $\hat f_{\vep|\sbX}(v|\bx)$, $F_{\vep|\sbX}(v|\bx)$
and $\hat \sigma(\bx)^2$. The purpose of this subsection is thus to
propose a way to estimate these factors and investigate their
asymptotic properties.

As pointed out by our referee, estimating $f_{\vep|\sbX}(0)$ is not a trivial task. The application of a rank test described in Chapter 3.5 of \cite{K:2005} is an alternative to avoid estimating $f_{\vep|\sbX}(0)$ in parametric quantile regression. However, it is a  challenging task to apply this technique to kernel smoothing quantile regression. For pointwise nonparametric inference, it may be possible to construct a test by adding weights (given by $h^{-1} K((\bx-\bX_i)/h)$, where $h$ is the bandwidth and $K$ is the kernel function) in the linear programing problem and therefore its dual can also be computed. However, a global shape test like the one investigated in this paper cannot be derived from the rank test. Hence, it seems inevitable to estimate the nuisance parameters and plug them into the test statistics.

Since we consider the additive error model \eqref{model},  the
conditional distribution function $F_{Y|\sbX}(\theta_0(\bx)|\bx)$
and the conditional density $f_{Y|\sbX}(\theta_0(\bx)|\bx)$ can be
replaced by $F_{\varepsilon|\sbX}(0|\bx)$ and
$f_{\vep|\sbX}(0|\bx)$, respectively, where $F_{\varepsilon|\sbX}$
and $f_{\varepsilon|\sbX}$ are the conditional distribution and
density functions of $\varepsilon$. Similarly, we have
$$
\sigma^2(\bx) = \E\big[\psi_\tau\big(Y-\theta_0(\bx)\big)^2 \big|\bX=\bx\big] = \E\big[\psi_\tau(\varepsilon)^2 \big|\bX=\bx\big]
$$
where $\vep$ may depend on $\bX$ due to heterogeneity. It should be
noted that the  kernel estimators for
$f_{\vep|\sbX}(0|\bx)$ and $f_{Y|\sbX}(\theta_0(\bx)|\bx)$ are asymptotically equivalent, but show different finite
sample behavior. We explore this issue further in the following
section.

Introducing the residuals $\hat \vep_i = Y_i-\hat \theta_n(\bX_i)$ we
propose to estimate $F_{\varepsilon|\sbX}$, $f_{\varepsilon|\sbX}$
and $\sigma^2(\bx)$ by
\begin{align}
\hat F_{\varepsilon|\sbX}(v|\bx) &= n^{-1}\sum_{i=1}^n G\left(\frac{v-\hat \vep_{i}}{h_0}\right) L_{\brh}(\bx-\bX_i)/\hat f_\sbX(\bx) \label{hat.F.eps.x},\\
\hat f_{\varepsilon|\sbX}(v|\bx) &= n^{-1}\sum_{i=1}^n g_{h_0}\left(v-\hat \vep_{i}\right) L_{\brh}(\bx-\bX_i)/\hat f_\sbX(\bx) \label{hat.f.eps.x},\\
\hat \sigma^2(\bx) &=  n^{-1}\sum_{i=1}^n \psi^2(\hat \vep_i) L_{\brh}(\bx-\bX_i)/\hat f_\sbX(\bx) \label{hat.sig},
\end{align}
where $\hat f_\sbX(\bx) = n^{-1} \sum_{i=1}^n L_{\brh}(\bx-\bX_i)$,
$G$ is a given continuously differentiable cumulative distribution
function and $g$ is its derivative. The construction of estimators in \eqref{hat.F.eps.x} and \eqref{hat.f.eps.x} follows from the estimator for general conditional distribution and density functions discussed in Chapter 5 and 6 of \cite{LR:2007}.
The same bandwidth $\brh$ is applied to the three estimators, but
the choice of $\brh$ will make the convergence rate of
\eqref{hat.sig} sub-optimal. More details on the choice of $\brh$ are given
in section \ref{subSec:bootQR} below. Nevertheless, the rate of convergence of \eqref{hat.sig} is of polynomial order in $n$. The theory developed in this
subsection can be generalized to the case of different bandwidth for
different direction without much difficulty.

The estimators \eqref{hat.F.eps.x} and \eqref{hat.f.eps.x} belong to
the family of residual-based estimators. The consistency of
residual-based density estimators for errors in a regression model
are explored in the literature in various settings. It is possible
to obtain an expression for the residual based kernel density
estimator as the sum of the estimator with the true residuals, the
partial sum of the true residuals and a term for the bias of the
nonparametrically estimated function, as shown in \cite{MN:2010},
among others. The residual based conditional kernel density case is
less considered in the literature. \cite{KN:2012} consider the
residual based kernel estimator for conditional distribution
function conditioning on a one-dimensional variable.

Below we give consistency results for the estimators defined in
\eqref{hat.F.eps.x}, \eqref{hat.f.eps.x} and \eqref{hat.sig}. The
proof can be found in the appendix.
\begin{lemma}\label{nuis.unifconv}
Under conditions (A1), (A3)-(A5), (B1)-(B3) in Section \ref{Sec:Ass}, we have
\begin{enumerate}
\item[1)] $\sup_{v \in I} \sup_{\sbx \in \mathcal D}\big|\hat
F_{\varepsilon|\sbX}(v|\bx)-F_{\varepsilon|\sbX}(v|\bx) \big|=
\CO_p\big(t_n\big)$,
\item[2)] $\sup_{v \in I} \sup_{\sbx
\in \mathcal D}\big|\hat
f_{\varepsilon|\sbX}(v|\bx)-f_{\varepsilon|\sbX}(v|\bx)\big|=\CO_p\big(t_n\big)$,
\item[3)] $\sup_{\sbx \in \mathcal D}\big|\hat
\sigma^2(\bx)-\sigma^2(\bx)\big| =
\CO_p\big(u_n\big)$,
\end{enumerate}
where $t_n = \CO\big\{h_0^{s'} + h^s + \brh^{s'} +
(n\brh^d)^{-1/2}\log n +(nh^d)^{-1/2}\log n\big\}=\CO(n^{-\lambda})$, and
$u_n= \CO\big\{h^s + \brh^{s'} + (n\brh^d)^{-1/2}\log n
+(nh^d)^{-1/2}\log n\big\}=\CO(n^{-\lambda_1})$ for some constants $\lambda, \lambda_1>0$.
\end{lemma}
The factor of $\log n$ shown in the convergence rate is the
price which we pay for the supnorm deviation. Since these estimators
uniformly converge in a polynomial rate in $n$, the asymptotic
distributions in Theorem \ref{TheoremQuantile} and
\ref{TheoremExpectile} do not change if we plug these estimators
into the formulae.


\section{Bootstrap confidence corridors}\label{Sec:Bootstrap}

\subsection{Asymptotic theory}\label{subSec:bootstrap}
In the case of the suitably normed maximum of independent standard normal variables, it
is shown in \cite{Hall:1979} that the speed of convergence in limit theorems of the form \eqref{max_dev}
is of order $1/\log n$, that is, the coverage error of the
asymptotic CC decays only logarithmically. This leads to
unsatisfactory finite sample performance of the asymptotic methods,
especially for small sample sizes and dimensions $d>1$. However, \cite{H:1991} suggests
that the use of a  bootstrap method, based on a proper way of
resampling, can increase the speed of shrinking of coverage error to
a polynomial rate of $n$. In this section we therefore propose a
specific bootstrap technique and construct a confidence corridor for
the objects to be analysed.

Given the residuals $\hat \vep_i = Y_i - \hat \theta_n(\bX_i)$, the
bootstrap observations $(\bX_i^*,\vep_{i}^*)$ are sampled from
\begin{align}
\hat f_{\vep,\sbX}(v,\bx) = \frac{1}{n}\sum_{i=1}^n g_{h_0}\left(\hat \vep_{i}-v\right) L_{\brh}(\bx-\bX_i), \label{fe_x}
\end{align}
where $g$ and $L$ are a kernel functions with bandwidths $h_0$,
$\brh$ satisfying assumptions (B1)-(B3). In particular, in our
simulation study, we choose $L$ to be a product Gaussian kernel. In
the following discussion $\P^*$ and $\E^*$ stand for the probability
and expectation conditional on the data
$(\bX_i,Y_i)$, $i=1,..., n$. 

We introduce the notation
$$
A_n^*(\bx) = \frac{1}{n}\sum_{i=1}^n K_h(\bx-\bX_i^*)\psi_\tau(\vep_i^*),
$$
and  define the so-called "one-step estimator" $\theta^*(\bx)$ from
the bootstrap sample by
\begin{align}
\hat \theta^*(\bx)-\hat \theta_n(\bx) =\hat S_{n,0,0}^{-1}(\bx) \left\{A_n^*(\bx)-\E^*[A_n^*(\bx)]\right\}, \label{theta.star}
\end{align}
where
\begin{align}\label{S.hat}
\hat S_{n,0,0}(\bx) =
 \left\{
  \begin{array}{ll}
    \hat f_{\vep|\sbX}\big(0|\bx\big) \hat f_\sbX(\bx), & \hbox{quantile case;} \\
    2 \big\{\tau-\hat F_{\vep|\sbX}\big(0|\bx\big)(2 \tau-1)\big\}\hat f_{\sbX}(\bx), & \hbox{expectile case.}
  \end{array}
\right.
\end{align}
note that $\E^*[\hat \theta^*(\bx)-\hat \theta_n(\bx)]=0$, so $\hat
\theta^*(\bx)$ is unbiased for $\hat \theta_n(\bx)$ under $\E^*$. As
a remark, we note that undersmoothing is applied in our procedure
for two reasons: first,  the theory we developed so far is based
on undersmoothing; secondly, it is suggested in \cite{Hall:1992}
that undersmoothing is more effective than oversmoothing given that
the goal is to achieve coverage accuracy.

Note that the bootstrap estimate \eqref{theta.star} is motivated by
the smoothed bootstrap procedure proposed in \cite{CK:2003}. In
contrast to these authors we make use of the leading term of the
Bahadur representation. \cite{MKY:2013} also use the leading term of
a Bahadur representation proposed in \cite{GS:2012} to construct
bootstrap samples. \cite{SRH:2012} propose a bootstrap for quantile
regression based on oversmoothing, which has the drawback that it
requires iterative estimation, and oversmoothing is in general less
effective in terms of coverage accuracy.

For the following discussion define
\begin{align}
Y_n^*(\bx) = \frac{1}{\sqrt{h^d \hat f_{\sbX}(\bx)} \sigma_{*}(\bx)} \int \int K\left(\frac{\bx-\bu}{h}\right)\psi_\tau\big(v\big)dZ_n^*(v,\bu)
\end{align}
as the bootstrap analogue of the process \eqref{Y}, where
\begin{align}
Z_n^*(y,\bu) = n^{1/2} \left\{F_n^*(v,\bu)-\hat
F(v,\bu)\right\}, \quad  \sigma_*(\bx)=\sqrt{\E^*\big[\psi_\tau(\vep_i^*)^2|\bx\big]} \label{sig.star}
\end{align}
and
$$
F_n^*(v,\bu) = \frac{1}{n} \sum_{i=1}^n \IF\left\{\vep_i^* \leq v, X_1^* \leq u_1, ..., X_d^* \leq u_d\right\}.
$$
The process $Y_n^*$ serves as an approximation of a standardized version of $\hat \theta_n^*-\hat \theta_n$, and similar to the previous sections the process
$Y_n^*$ is approximated by a stationary Gaussian field $Y_{n,5}^*$ under
$\P^*$ with probability one, that is,
\begin{align*}
  Y_{5,n}^*(\bx) = \frac{1}{\sqrt{h^d}} \int
K\left(\frac{\bx-\bu}{h}\right)dW^*(\bu).
\end{align*}
Finally, $\sup_{\sbx \in \mathcal D}\big|Y_{5,n}^*(\bx)\big|$ is
asymptotically Gumbel distributed conditional on samples.

\begin{theorem}\label{TheoremBoot}
Suppose that assumptions (A1)-(A6), (C1) in Section \ref{Sec:Ass} hold, and  
$\mbox{vol}(\mathcal D)=1$, let
\begin{align*}
r_n^*(\bx) = \sqrt{\frac{nh^d}{\hat f_\sbX(\bx) \sigma_*^2(\bx)}} \hat S_{n,0,0}(\bx),
\end{align*}
where $\hat S_{n,0,0}(\bx)$ is defined in \eqref{S.hat} and
$\sigma_*^2(\bx)$ is defined in \eqref{sig.star}. Then
\begin{align}
\P^*\left\{(2d \kappa \log n)^{1/2} \left(\sup_{\sbx \in \mathcal D}\big[r_n^*(\bx)|\hat \theta^*(\bx)-\hat \theta_n(\bx)|\big]/\|K\|_2 -d_n\right) < a \right\} \rightarrow \exp\big\{-2\exp(-a)\big\}, \quad \mbox{a.s.} \label{boot.theorem}
\end{align}
as $n \rightarrow \infty$ for the local constant quantile regression
estimate. If (A1)-(A6) and (EC1) hold with a constant $b \geq 4$
satisfying
$$
n^{-\frac{1}{6}+\frac{4}{b^2}-\frac{1}{b}}h^{-\frac{d}{2}-\frac{6d}{b}} = \CO(n^{-\nu}), \quad \nu >0,
$$
then \eqref{boot.theorem} also holds for expectile regression with
corresponding $\sigma_*^2(\bx)$.
\end{theorem}
The proof can be found in the appendix. The following lemma suggests
that we can replace $\sigma_\ast(\bx)$ in the
limiting theorem by $\hat \sigma(\bx)$. 
\begin{lemma}\label{sigstar.consistent}
If assumptions  (B1)-(B3), and (EC1) in Section \ref{Sec:Ass} are satisfied with $b>2(2s'+d+1)/(2s'+3)$, then
$$\|\sigma_*^2(\bx)-\hat \sigma^2 (\bx)\|=\Co_p^*\big((\log
n)^{-1/2}\big), \quad a.s.$$
\end{lemma}

The following corollary is a consequence of Theorem \ref{TheoremBoot}.
\begin{corollary}
  Under the same conditions as stated in Theorem \ref{TheoremBoot}, the (asymptotic) bootstrap confidence set of level $1-\alpha$ is given by
  \begin{align}
  \left\{\theta:\sup_{\sbx\in\mathcal D}\left|\frac{\hat S_{n,0,0}(\bx)}{\sqrt{\hat f_\sbX(\bx) \hat \sigma^2(\bx)}} \big[\hat \theta_n(\bx)-\theta(\bx)\big]\right| \leq \xi_\alpha^* \right\}, \label{bootdist.er}
\end{align}
  where $\xi_\alpha^*$ satisfies
  \begin{align}
 \lim_{n\rightarrow\infty}
\P^*\left(\sup_{\sbx\in\mathcal D}\left|\frac{\hat S_{n,0,0}(\bx)}{\sqrt{\hat f_\sbX(\bx) \hat \sigma^2(\bx)}} \big[\hat \theta^*(\bx)-\hat \theta_n(\bx)\big]\right|\leq \xi_\alpha^* \right)=1-\alpha, \quad a.s.
  \end{align}
where $\hat S_{n,0,0}$ is
defined in \eqref{S.hat}.
\end{corollary}

Note that it does not create much difference to standardize the
$\hat \theta_n(\bx)-\theta_0(\bx)$ in \eqref{boot.theorem} with
$\hat f_\sbX$ and $\hat \sigma^2(\bx)$ constructed from original
samples or $\hat f_\sbX$ and $\hat \sigma^2(\bx)$ from the bootstrap samples. 
The simulation results of \cite{CK:2003} show that the two ways of
standardization give similar coverage probabilities for confidence
corridors of kernel ML estimators.

\subsection{Implementation}\label{subSec:bootQR}
In this section, we discuss issues related to the implementation of
the bootstrap for quantile regression.

Note that the \emph{width} of the CC is determined by the variance and the \emph{location} is affected by the bias of the quantile function estimator, and both depend on the bandwidth used for estimation. Hence, the choice of bandwidth needs to balance the bias (location) and the variance (size). It is chosen such that the bias is only just negligible after normalization, that is, slightly smaller than the $L^2$-optimal bandwidth. Therefore, it is enough to take an undersmoothed $h
=\CO(n^{-1/(2s+d)-\delta})$, given that $s>d$ and $\delta>0$, where $s$ is the order of H\"older continuity of the function $\theta_0$ and $\delta$ is the degree of undersmoothing. We may use the 
methods proposed by \cite{YJ:1998} for nonparametric quantile regression to choose the bandwidth before undersmoothing, namely
\begin{align}
	h_{\tau,j} = h_{1,j}\{\tau(1-\tau)/\phi(\Phi^{-1}(\tau))^2\}^{1/5}, \quad j=1,2, \label{yujones}	
\end{align}
where $h_{1,j}$ are chosen by common methods like the rule-of-thumb or cross-validation for mean regression or density estimation and $\Phi$ is the CDF of the standard Gaussian distribution. In our simulation study, we select $h_{1,j}$ in \eqref{yujones} by the rule-of-thumb, implemented with the \texttt{np} package in \texttt{R}. In our application analysis, $h_{1,j}$ in \eqref{yujones} are chosen by the cross-validated bandwidth for the conditional distribution smoother of $Y$ given $\bX$, implemented with the \texttt{np} package in \texttt{R}. This package is based on the paper of \cite{LLR:2013}.

For expectile regression, we use the rule-of-thumb bandwidth for the conditional distribution smoother of $Y$ given $\bX$, chosen with the \texttt{np} package in \texttt{R}. 

The choice of $h_0$ and $\brh$ for estimating the scaling factors in Section \ref{subSec:scaling} should minimize the convergence rate
of these residual based estimators. Hence, observing that the terms
 related to $h_0$ and $\brh$ are similar to those in usual
$(d+1)$-dimensional density estimators, it is reasonable to choose $h_0 \sim \brh
\sim n^{-1/(5+d)}$, given that $L$, $g$ are second order kernels. We
choose the rule-of-thumb bandwidths for conditional densities with
the \texttt{R} package \texttt{np} in our simulation and application
studies.

The one-step estimator for quantile regression defined in
\eqref{theta.star} depends sensitively on the estimator of $\hat
S_{n,0,0}(\bx)$. Unlike in the expectile case, the function
$\psi(\cdot)$ in the quantile case is bounded, and, as a result, the
bootstrapped density based on \eqref{bootdist.er} is very easily
influenced by the factor $\hat S_{n,0,0}(\bx)$; in particular, $\hat
f_{\vep|\sbX}(0|\bx)$. As pointed out by \cite{FHH:2011}, the
residual of quantile regression tends to be less dispersed than the
model error; thus $\hat f_{\vep|\sbX}(0|\bx)$ tends to over-estimate
the true $f_{\vep|\sbX}(0|\bx)$ for each $\bx$.


The way of getting around this problem is based on the following
observation: An additive error model implies the equality
$f_{Y|\sbX}\big\{v+\theta_0(\bx)|\bx\big\} =
f_{\vep|\sbX}\big(v|\bx\big)$, but this property does not hold for the kernel estimators
\begin{align}
\hat f_{\varepsilon|\sbX}(0|\bx) &= n^{-1}\sum_{i=1}^n g_{h_0}\left(\hat \vep_{i}\right) L_{\brh}(\bx-\bX_i)/\hat f_\sbX(\bx), \\
\hat f_{Y|\sbX}(\hat \theta_n(\bx)|\bx) &= n^{-1}\sum_{i=1}^n g_{h_1}\left(Y_i-\hat \theta_n(\bx)\right) L_{\tilde h}(\bx-\bX_i)/\hat f_\sbX(\bx), \label{fy_x}
\end{align}
 of the conditional density functions.
In general $\hat f_{\varepsilon|\sbX}(0|\bx) \neq \hat
f_{Y|\sbX}(\hat \theta_n(\bx)|\bx)$ in $\bx$ although both estimates
are asymptotically equivalent. In applications the two estimators
can differ substantially due to the bandwidth selection because for
data-driven bandwidths we usually have $h_0 \neq h_1$. For example, if a 
common method for bandwidth selection such as a rule-of-thumb is used,
$h_1$ will tend to be larger than $h_0$ since the sample variance of
$Y_i$ tends to be larger than that of $\hat \vep_i$. Given that the
same kernels are applied, it happens often that $\hat
f_{Y|\sbX}(\hat \theta_n(\bx)|\bx) > f_{Y|\sbX}(\theta_0(\bx)|\bx)$,
even if $\hat \theta_n(\bx)$ is usually very close to
$\theta_0(\bx)$. To correct such abnormality, we are motivated to
set $h_1=h_0$ which is the rule-of-thumb bandwidth of $\hat
f_{\vep|\sbx}(v|\bx)$ in \eqref{fy_x}. As the result, it leads to a
more rough estimate for $\hat f_{Y|\sbX}(\hat \theta_n(\bx)|\bx)$.

In order to exploit the roughness of $\hat f_{Y|\sbX}(\hat
\theta_n(\bx)|\bx)$ while making the CC as narrow as possible, we
develop a trick depending on
\begin{align}
\frac{\hat f_{Y|\sbX}\big\{\hat \theta_n(\bx)|\bx\big\}}{\hat f_{\vep|\sbX}(0|\bx)} = \frac{h_0}{h_1}\frac{\sum_{i=1}^n g_{h_1}\left(\big\{Y_i-\hat \theta_n(\bx)\big\}/h_1\right) L_{\tilde h}(\bx-\bX_i)}{\sum_{i=1}^n g_{h_0}\left(\hat \vep_{i}/h_0\right) L_{\brh}(\bx-\bX_i)}. \label{densratio}
\end{align}
As $n \to \infty$, \eqref{densratio} converges to 1. If we impose
$h_0=h_1$, as the multiple $h_0/h_1$ vanishes, \eqref{densratio}
captures the deviation of the two estimators without the difference
of the bandwidth in the way. In particular,
the bandwidth $h_0=h_1$ is selected as the rule-of-thumb bandwidth
for $\hat f_{\vep|\sbX}(y|\bx)$. This makes $\hat
f_{\vep|\sbX}(y|\bx)$ larger and thus leads to a narrower CC, as
will be more clear below.

We propose the alternative bootstrap confidence corridor for
quantile estimator:
  $$
  \left\{\theta:\sup_{\sbx\in\mathcal D}\big|\sqrt{\hat f_\sbX(\bx)}\hat f_{Y|\sbX}\big\{\hat \theta_n(\bx)|\bx\big\}\big[\hat \theta_n(\bx)-\theta(\bx)\big]\big| \leq \xi_\alpha^\dag \right\},
  $$
where $\xi_\alpha^\dag$ satisfies
\begin{align}
  \P^*\left(\sup_{\sbx\in\mathcal D}\left|\hat f_\sbX(\bx)^{-1/2} \frac{\hat f_{Y|\sbX}\big\{\hat \theta_n(\bx)|\bx\big\}}{\hat f_{\vep|\sbX}(0|\bx)}\big[A_n^*(\bx)-\E^*A_n^*(\bx)\big]\right|\leq \xi_\alpha^\dag \right)=1-\alpha.\label{qr.boot.cc}
  \end{align}
Note that the probability on the left-hand side of
\eqref{qr.boot.cc} can again be approximated by a Gumbel
distribution function asymptotically, which follows by Theorem
\ref{TheoremBoot}.

\section{A simulation study}\label{Sec:simu}
In this section we investigate the methods described in the previous
sections by means of a simulation study.  We construct confidence corridors
for quantiles and expectiles for different levels $\tau$ and use the
quartic (product) kernel. The performance of our methods is compared to the performance of the method proposed by \cite{HH:12} at the end of this section.
 For the confidence based on asymptotic
distribution theory, we use the rule of thumb bandwidth chosen from
the \texttt{R} package \texttt{np}, and then rescale it as described
in \cite{YJ:1998}, finally multiply it by $n^{-0.05}$ for
undersmoothing. The sample sizes are given by $n=100, 300$ and 500, so
the undersmoothing multiples are $0.794$, $0.752$ and $0.733$
respectively. We take $20 \times 20$ equally distant grids in $[0.1,0.9]^2$ and estimate quantile or expectile functions pointwisely on this set of grids. In the quantile regression bootstrap CC, the bandwidth
$h_1$ used for estimating $\hat f_{Y|\sbX}(y|\bx)$ is chosen to be
the rule-of-thumb bandwidth of $\hat f_{\vep|\sbX}(0|\bx)$ and
multiplied by a multiple 1.5. This would give slightly wider CCs.

\begin{table}[h!]
        \begin{center}\fontsize{11}{12} \selectfont
        \begin{tabular}{llllllll}
\hline\hline
& &\multicolumn{3}{c}{\textbf{Homogeneous}}  &     \multicolumn{3}{c}{\textbf{Heterogeneous}}  \\
\textbf{Method} &\boldmath $n$ & \boldmath $\tau=0.5$ & \boldmath $\tau=0.2$  & \boldmath $\tau=0.8$ & \boldmath $\tau=0.5$ & \boldmath $\tau=0.2$  & \boldmath $\tau=0.8$\\
\hline
\multicolumn{8}{c}{\boldmath $\sigma_0 = 0.2$}    \\[0.7ex]
          &100  &.000(0.366) &.109(0.720) &.104(0.718) &.000(0.403) &.120(0.739) &.122(0.744)\\
          & 300 & .000(0.304) &.130(0.518) &.133(0.519) &.002(0.349) &.136(0.535) &.153(0.537)\\
          &500  & .000(0.262) &.117(0.437) &.142(0.437) &.008(0.296) &.156(0.450)&.138(0.450)\\
 \multicolumn{8}{c}{\boldmath $\sigma_0 = 0.5$}    \\[0.7ex]
         & 100  & .070(0.890) & .269(1.155) & .281(1.155)  &.078(0.932) &.300(1.193) &.302(1.192) \\
 Asympt. &300 & .276(0.735)&.369(0.837) &.361(0.835) &.325(0.782) &.380(0.876) &.394(0.877) \\
         &500  &.364(0.636) &.392(0.711) &.412(0.712) &.381(0.669)&.418(0.743)&.417(0.742)\\
 \multicolumn{8}{c}{\boldmath $\sigma_0 = 0.7$}    \\[0.7ex]
         &100   &.160(1.260) &.381(1.522) &.373(1.519) &.155(1.295) &.364(1.561) &.373(1.566) \\
         & 300&.438(1.026) &.450(1.109) &.448(1.110) &.481(1.073) &.457(1.155) &.472(1.152) \\
         &500    &.533(0.888) &.470(0.950) &.480(0.949) &.564(0.924)&.490(0.984)&.502(0.986)\\[0.7ex]
                \hline
                \multicolumn{8}{c}{\boldmath $\sigma_0 = 0.2$}    \\[0.7ex]
        & 100&.325(0.676) &.784(0.954) &.783(0.954) &.409(0.717) &.779(0.983) &.778(0.985)\\
        &300 &.442(0.457) &.896(0.609) &.894(0.610) &.580(0.504) &.929(0.650) &.922(0.649)\\
        &500 &.743(0.411) &.922(0.502) &.921(0.502) &.839(0.451) &.950(0.535) &.952(0.536)\\

 \multicolumn{8}{c}{\boldmath $\sigma_0 = 0.5$}    \\[0.7ex]
          &100  &.929(1.341) &.804(1.591) &.818(1.589) &.938(1.387) &.799(1.645) &.773(1.640) \\
   Bootst.&300 &.950(0.920)&.918(1.093) &.923(1.091) &.958(0.973) &.919(1.155) &.923(1.153) \\
          & 500 &.988(0.861) &.968(0.943) &.962(0.942) &.990(0.902) &.962(0.986) &.969(0.987) \\

 \multicolumn{8}{c}{\boldmath $\sigma_0 = 0.7$}    \\[0.7ex]
        &100 &.976(1.811) &.817(2.112) &.808(2.116) &.981(1.866) &.826(2.178) &.809(2.176) \\
        &300 &.986(1.253) &.919(1.478) &.934(1.474) &.983(1.308) &.930(1.537) &.920(1.535) \\
        &500 &.996(1.181) &.973(1.280) &.968(1.278) &.997(1.225) &.969(1.325) &.962(1.325) \\
\hline\hline
        \end{tabular}
        \end{center}
\caption{{\it Nonparametric quantile model coverage probabilities.
The nominal coverage is $95\%$. The number in the parentheses is the
volume of the confidence corridor. The asymptotic method corresponds to the asymptotic quantile regression CC and bootstrap method corresponds to quantile regression bootstrap CC.}
}\label{qrratio.asym}
        \end{table}

\begin{table}[h!]
        \begin{center}\fontsize{11}{12} \selectfont
        \begin{tabular}{llllllll}
\hline\hline
& &\multicolumn{3}{c}{\textbf{Homogeneous}}  &     \multicolumn{3}{c}{\textbf{Heterogeneous}}  \\
\textbf{Method} &\boldmath $n$ & \boldmath $\tau=0.5$ & \boldmath $\tau=0.2$  & \boldmath $\tau=0.8$ & \boldmath $\tau=0.5$ & \boldmath $\tau=0.2$  & \boldmath $\tau=0.8$\\
\hline
\multicolumn{8}{c}{\boldmath $\sigma_0 = 0.2$}    \\[0.7ex]
         &100 & .000(0.428) &.000(0.333) &.000(0.333) &.000(0.463) &.000(0.362) &.000(0.361) \\
         &300 & .049(0.341) &.000(0.273) &.000(0.273) &.079(0.389) &.001(0.316) &.002(0.316) \\
         &500 & .168(0.297) &.000(0.243) &.000(0.243) &.238(0.336) &.003(0.278) &.002(0.278) \\
\multicolumn{8}{c}{\boldmath $\sigma_0 = 0.5$}    \\[0.7ex]
         &100 & .007(0.953) &.000(0.776)&.000(0.781) &.007(0.997)&.000(0.818)&.000(0.818) \\
 Asympt. &300 &.341(0.814)&.019(0.708)&.017(0.709)&.355(0.862)&.017(0.755)&.018(0.754)\\
         &500 &.647(0.721)&.067(0.645)&.065(0.647) &.654(0.759)&.061(0.684)&.068(0.684)\\
\multicolumn{8}{c}{\boldmath $\sigma_0 = 0.7$}    \\[0.7ex]
        &100 & .012(1.324) &.000(1.107)&.000(1.107) &.010(1.367)&.000(1.145)&.000(1.145) \\
        &300 & .445(1.134)&.021(1.013)&.013(1.016)&.445(1.182)&.017(1.062)&.016(1.060)\\
        & 500 & .730(1.006) &.062(0.928)&.078(0.929) &.728(1.045)&.068(0.966)&.066(0.968)\\[0.7ex]
    \hline
    \multicolumn{8}{c}{\boldmath $\sigma_0 = 0.2$}    \\[0.7ex]
          & 100 &.686(2.191) &.781(2.608) &.787(2.546) &.706(2.513) &.810(2.986) &.801(2.943) \\
        & 300 &.762(0.584) &.860(0.716) &.876(0.722) &.788(0.654) &.877(0.807) &.887(0.805) \\
        & 500 &.771(0.430) &.870(0.533) &.875(0.531) &.825(0.516) &.907(0.609) &.904(0.615) \\
\multicolumn{8}{c}{\boldmath $\sigma_0 = 0.5$}    \\[0.7ex]
        & 100 &.886(5.666) &.906(6.425)&.915(6.722) &.899(5.882)&.927(6.667)&.913(6.571) \\
 Bootst. &300&.956(1.508)&.958(1.847) &.967(1.913) &.965(1.512)&.962(1.866)&.969(1.877)\\
        & 500 &.968(1.063)&.972(1.322) &.972(1.332) &.972(1.115)&.971(1.397)&.974(1.391)\\
\multicolumn{8}{c}{\boldmath $\sigma_0 = 0.7$}    \\[0.7ex]
        & 100 &.913(7.629) &.922(8.846) &.935(8.643) &.929(8.039) &.935(9.057) &.932(9.152)     \\
        & 300 &.969(2.095) &.969(2.589) &.971(2.612) &.974(2.061) &.972(2.566) &.979(2.604)\\
          & 500 &.978(1.525) &.976(1.881) &.967(1.937) &.981(1.654) &.978(1.979) &.974(2.089)\\
\hline\hline
        \end{tabular}
        \end{center}
\caption{{\it Nonparametric expectile model coverage probability. The nominal
coverage is $95\%$. The number in the parentheses is the volume of
the confidence corridor.  The
asymptotic method corresponds to the asymptotic expectile regression
CC and bootstrap method corresponds to expectile regression
bootstrap CC.}}\label{erratio.asym}
        \end{table}

The data are generated from the normal regression model
\begin{align}
Y_i = f(X_{1,i},X_{2,i}) + \sigma(X_{1,i},X_{2,i}) \varepsilon_i, \quad i=1,\ldots,n \label{sim.model}
\end{align}
where the independent variables $(X_1,X_2)$ follow a joint uniform
distribution taking values on $[0,1]^2$, $\cov(X_1,X_2)=0.2876$,
$f(X_1,X_2) = \sin(2 \pi X_1)+X_2$, and $\varepsilon_i$ are
independent standard Gaussian random variables. For both quantile
and expectile, we look at three quantiles of the distribution, namely
$\tau=0.2,0.5,0.8$. The set of grid point is $G \times G$ where $G$ is the set of 20 equidistant grids on univariate interval $[0.1,0.9]$. Thus, the grid size is $|G \times G|=400$. 

In the homogeneous model, we take $\sigma(X_1,X_2)=\sigma_0$, for
$\sigma_0 = 0.2, 0.5, 0.7$. In the heterogeneous model, we take
$\sigma(X_1,X_2) = \sigma_0 + 0.8 X_1(1-X_1)X_2(1-X_2)$. 2000
simulation runs are carried out to estimate the coverage probability.

The upper part of Table \ref{qrratio.asym} shows the coverage probability of the asymptotic CC
for nonparametric quantile regression functions. It can be
immediately seen that the asymptotic CC performs very poorly,
especially when $n$ is small. A comparison of the results with
those of one-dimensional asymptotic simultaneous confidence bands derived in
\cite{CK:2003} or \cite{FL:2013}, shows that the accuracy in the two-dimensional case is much worse. Much to our surprise, the asymptotic CC
performs better in the case of $\tau=0.2,0.8$ than in the case of
$\tau=0.5$. On the other hand, it is perhaps not so amazing to see
that asymptotic CCs behave similarly under both homogeneous and
heterogeneous models. As a final remark about the asymptotic CC we mention
that it is highly sensitive with respect to $\sigma_0$. Increasing values of $\sigma_0$ yields
larger CC, and this may lead to greater coverage probability.

The lower part of Table \ref{qrratio.asym} shows that the bootstrap CCs for nonparametric quantile regression functions yield a remarkable improvement in comparison to the
 asymptotic CC. For the bootstrap CC, the coverage
probabilities are in general close to the nominal coverage of 95\%.
The bootstrap CCs are usually wider, and getting narrower when $n$
increases. Such phenomenon can also be found in the simulation study
of \cite{CK:2003}. Bootstrap CCs are less sensitive than asymptotic
CCs with respect to the choice $\sigma_0$, which is also considered
as an advantage. Finally, we note that the performance of bootstrap
CCs does not depend on which variance specification is used too.

The upper part of Table \ref{erratio.asym} shows the coverage probabiltiy of the CC
for nonparametric expectile regression functions. The results are similar to the case of quantile regression. The asymptotic CCs do \emph{not} give accurate coverage probabilities. For example in some cases like $\tau=0.2$ and $\sigma_0=0.2$, not a single simulation in the 2000 iterations yields a case where surface is completely covered by the asymptotic CC. 

The lower part of Table \ref{erratio.asym} shows that bootstrap CCs for expectile regression give more accurate approximates to the nominal coverage than the asymptotic CCs. One can see in the parenthesis that the volumes of the bootstrap CCs are significantly larger than those of the asymptotic CCs, especially for small $n$. 

\begin{table}[h!]
        \begin{center}\fontsize{11}{12} \selectfont
        \begin{tabular}{lllllll}
\hline\hline
 &\multicolumn{3}{c}{\textbf{Homogeneous}}  &     \multicolumn{3}{c}{\textbf{Heterogeneous}}  \\
 \boldmath $n$ & \boldmath $\xi=0.005$ & \boldmath $\xi=0.05$  &\boldmath $\xi=0.1$ & \boldmath $\xi=0.005$ & \boldmath $\xi=0.05$  &\boldmath $\xi=0.1$\\
\hline
\multicolumn{7}{c}{\boldmath $\sigma_0 = 0.2$}    \\[0.7ex]
         100 &.693(3.027) &.529(1.740) &.319(1.040) &.680(3.452) &.546(2.051) &.332(1.224) \\
         300 &.891(0.580) &.748(0.365) &.642(0.323) &.907(0.667) &.798(0.414) &.698(0.364) \\
         500 &.886(0.335) &.770(0.265) &.678(0.244) &.896(0.379) &.789(0.298) &.699(0.274) \\
\multicolumn{7}{c}{\boldmath $\sigma_0 = 0.5$}    \\[0.7ex]
         100 &.720(7.264) &.611(4.489) &.394(2.686) &.729(7.594) &.616(4.676) &.414(2.829)  \\
         300 &.945(1.423) &.849(0.859) &.755(0.746) &.940(1.511) &.854(0.912) &.760(0.791) \\
         500 &.944(0.795) &.846(0.600) &.750(0.548) &.937(0.833) &.839(0.632) &.751(0.577) \\
\multicolumn{7}{c}{\boldmath $\sigma_0 = 0.7$}    \\[0.7ex]
        100 &.730(10.183) &.634(6.411) &.430(3.853) &.752(10.657)&.658(6.577) &.441(3.923) \\
        300 &.936(1.995) &.854(1.197)  &.751(1.037) &.951(2.091) &.875(1.256) &.772(1.086) \\
        500 &.933(1.098) &.854(0.831)  &.774(0.758) &.938(1.145) &.853(0.865) &.770(0.789)\\[0.7ex]
\hline\hline
        \end{tabular}
        \end{center}
\caption{{\it Proportion in 2000 iteration that the coverage of $\geq 95\%$ grid points for nonparametric mean model, using the bootstrap method of \cite{HH:12}. The nominal
coverage is $95\%$. The number in the parentheses is the volume of
the confidence corridor.}}\label{erratio.hh}
        \end{table}

Table \ref{erratio.hh} presents the proportion in the 2000 iterations which covers 95\% of the 400 grid points, using the bootstrap method proposed in \cite{HH:12}(abbreviated as HH) for nonparametric mean regression at $d=2$. HH derived an expansion for the
bootstrap bias and established a somewhat different way to construct
confidence bands without the use of extreme value theory. It is worth noting that their
bands are uniform with respect to a fixed but unspecified portion of $(1-\xi)\cdot100\%$
(smaller than $100\%$) of grid points, while in our approach the uniformity is achieved on
the whole set of grids.
 
The simulation model is \eqref{sim.model} with the same homogeneous and heterogeneous variance specifications as before. We choose three levels of $\xi = 0.005, 0.05$ and $0.1$. It is suggested in HH that $\xi=0.1$ is usually sufficient in univariate nonparametric mean regression $d=1$. Note that $\xi=0.005$ corresponds to the second smallest pointwise quantile $\hat\beta(\bx,0.05)$ in the notation of HH, given that our grid size is 400. This is close to the uniform CC in our sense. The simulation model associated with the Table \ref{erratio.hh} is the same with that of the case $\tau=0.5$ in the bootstrap part of Table \ref{qrratio.asym} and Table \ref{erratio.asym}, because in case of the normal distribution the median equals the mean and $\tau=0.5$ expectile is exactly the mean. However, one should be aware that our coverage probabilities are more stringent because we check the coverage at every point in the set of grids, rather than only 95\% of the points (we refer it as \emph{complete coverage}). Hence, the complete coverage probability of HH will be lower than the proportion of 95\% coverage shown in Table \ref{erratio.hh}. The proportion of 95\% coverage should therefore be viewed as an upper bound for the complete coverage.

We summarize our findings as follows. Firstly the proportion of 95\% coverage in general present similar patterns as shown in Table \ref{qrratio.asym} and \ref{erratio.asym}. The coverage improves when $n$ and $\sigma_0$ get larger, and the volume of the band decreases as $n$ increases and increases when $\sigma_0$ increases. The homogeneous and heterogeneous model yield similar performance. Comparing with the univariate result in HH, it is found that the proportion of coverage tends to perform worse than that in HH under the same sample size. This is due to the curse of dimensionality, the estimation of a bivariate function is less accurate than that of an univariate function. As the result, a more conservative $\xi$ has to be applied. If we compare Table \ref{erratio.hh} to the bootstrap part of \ref{qrratio.asym} with $\tau=0.5$, it can be seen that our complete coverage probabilities are comparable to the proportion of 95\% coverage at the case $\xi=0.005$, though in the case of $\sigma_0=0.2$ our CC does not perform very well. However, the volumes of our CC are much less than that of HH in the cases of small $n$ and moderate and large $\sigma_0$. This suggests that our CC is more efficient. Finally, the proportion of 95\% coverage at $\xi=0.005$ in Table \ref{erratio.hh} is similar to the complete coverage probability in bootstrap part of \ref{erratio.asym} with $\tau=0.5$, but when sample size is small, the volume of our CC is smaller. 
 
\section{Application: a treatment effect study}\label{Sec:apply}
The classical application of the proposed method consists in testing the
hypothetical functional form of the regression function. 
Nevertheless, the proposed method can also be applied to test for a
quantile treatment effect \citep[see][]{K:2005} or to test for
conditional stochastic dominance (CSD) as investigated in
\cite{DE:2013}. In this section we shall apply the new method to
test these hypotheses for data collected from a real government
intervention.

The estimation of the quantile treatment effect (QTE) recovers the
heterogeneous impact of intervention on various points of the
response distribution. To define QTE, given vector-valued exogenous
variables $\bX \in \bm{\mathcal{X}}$ where $\bm{\mathcal{X}} \subset
\R^d$, suppose $Y_0$ and $Y_1$ are response variables associated
with the control group and treatment group, and let $F_{0|\sbX}$ and
$F_{1|\sbX}$ be the conditional distribution for $Y_0$ and $Y_1$,
the QTE at level $\tau$ is defined by
\begin{align}
\Delta_\tau(\bx) \defeq Q_{1|\sbX}(\tau|\bx)-Q_{0|\sbX}(\tau|\bx), \quad \bx \in \bm{\mathcal{X}}, \label{def.qte}
\end{align}
where $Q_{0|\sbX}(y|\bx)$ and $Q_{1|\sbX}(y|\bx)$ are the
conditional quantile of $Y_0$ given $\bX$ and $Y_1$ given $\bX,$
respectively. This definition corresponds to the idea of horizontal
distance between the treatment and control distribution functions
appearing in \cite{Dok:1974} and \cite{Leh:1974}.

A related concept in measuring the efficiency of a treatment is the
so called "conditional stochastic dominance". $Y_1$ conditionally
stochastically dominates $Y_0$ if
\begin{align}
F_{1|\sbX}(y|\bx) \leq F_{0|\sbX}(y|\bx) \quad \mbox{a.s. for all }(y,\bx) \in (\mathcal Y, \bm{\mathcal{X}}), \label{csd}
\end{align}
where $\mathcal Y$, $\bm{\mathcal{X}}$ are domains of $Y$ and $\bX$.
For example, if $Y_0$ and $Y_1$ stand for the income of two groups
of people $G_0$ and $G_1$, \eqref{csd} means that the distribution
of $Y_1$ lies on the right of that of $Y_0$, which is equivalent to
saying that at a given $0<\tau<1$, the $\tau$-quantile of $Y_1$ is
greater than that of $Y_0$. Hence, we could replace the testing
problem \eqref{csd} by
\begin{align}
Q_{1|\sbX}(\tau|\bx) \geq Q_{0|\sbX}(\tau|\bx) \quad \mbox{ for all } 0<\tau<1 \mbox{ and } \bx \in \bm{\mathcal{X}}. \label{csd_quantile}
\end{align}
Comparing \eqref{csd_quantile} and \eqref{def.qte}, one would find
that \eqref{csd_quantile} is just a uniform version of the test
$\Delta_\tau(\bx) \geq 0$ over $0<\tau<1$.

The method that we introduced in this paper is suitable for testing
a hypothesis like $\Delta_\tau(\bx) = 0$ where $\Delta_\tau(\bx)$ is
defined in \eqref{def.qte}. One can construct CCs for
$Q_{1|\sbX}(\tau|\bx)$ and $Q_{0|\sbX}(\tau|\bx)$ respectively, and
then check if there is overlap between the two confidence regions.
One can also extend this idea to test \eqref{csd_quantile} by
building CCs for several selected levels $\tau$.

We use our method to test the effectiveness of the National
Supported Work (NSW) demonstration program, which was a randomized,
temporary employment program initiated in 1975 with the  goal to
provide work experience for individuals who face economic and social
problems prior to entering the program. The data have been widely
applied to examine techniques which estimate the treatment
effect in a nonexperimental setting. In a pioneer study,
\cite{lalonde:1986} compares the treatment effect estimated from the
experimental NSW data with that implied by nonexperimental
techniques. \cite{DW:1999} analyse a subset of Lalonde's data and
propose a new estimation procedure for nonexperimental treatment
effect giving more accurate estimates than Lalonde's estimates. The
paper that is most related to our study is \cite{DE:2013}. These authors
propose a test for hypothesis \eqref{csd} and apply it to
Lalonde's data, in which they choose "age" as the only conditional
covariate and the response variable being the increment of earnings
from 1975 to 1978. They cannot reject the null hypothesis of
nonnegative treatment effect on the earnings growth.
\begin{figure}[!h]
\centering
 \includegraphics[width=7.8cm, height = 7cm]{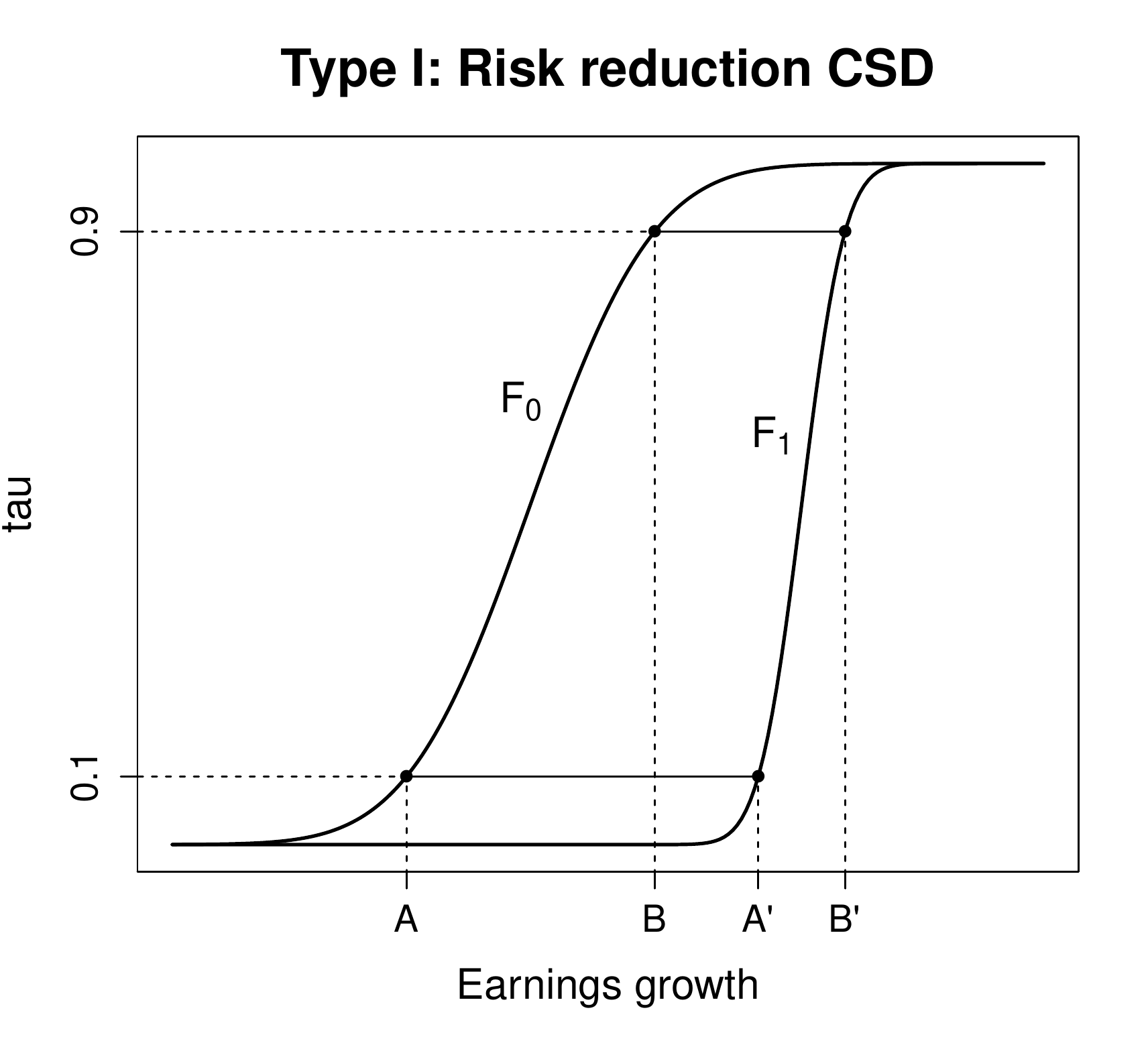}
 \includegraphics[width=7.8cm, height = 7cm]{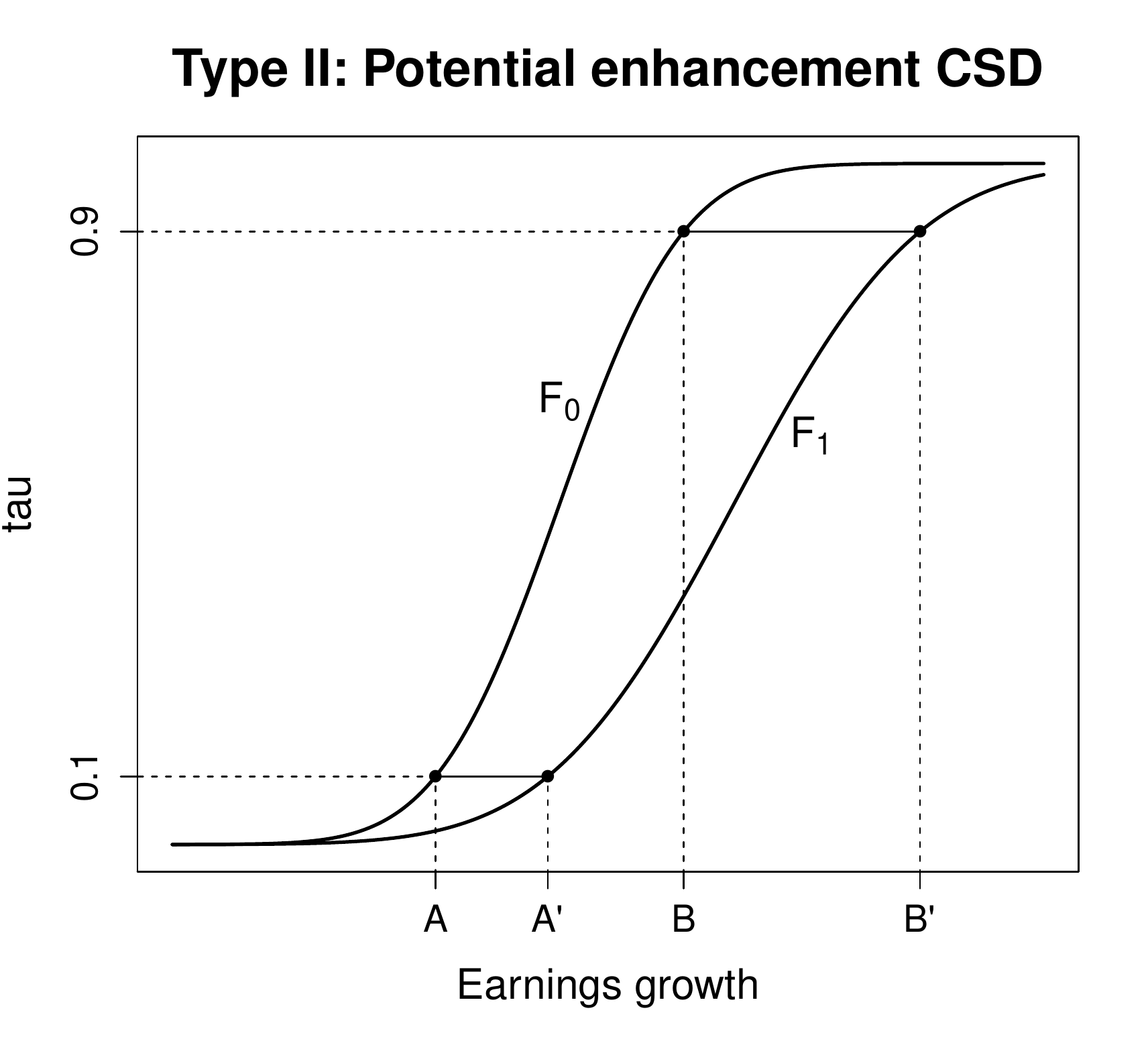}
\caption{The illustrations for the two possible types of stochastic dominance. In the left figure, the 0.1 quantile improves (downside risk reduction) more dramatically than the 0.9 quantile (upside potential increase), as the distance between $A$ and
$A'$ is greater than that between $B$ and $B'$. For the right picture the interpretation is just
the opposite.}\label{illus.sd}
\end{figure}

The previous literature, however, has not addressed an important
question. We shall depict this question by two pictures. In Figure
\ref{illus.sd}, it is obvious that $Y_1$ stochastically dominates
$Y_0$ in both pictures, but significant differences can be seen
between the two scenarios. For the left one, the 0.1 quantile improves more
dramatically than the 0.9 quantile, as the distance between $A$ and
$A'$ is greater than that between $B$ and $B'$. In usual words, the
gain of the 90\% lower bound of the earnings growth is more than
that of the 90\% upper bound of the earnings growth after the
treatment. "90\% lower bound of the earnings growth" means the
probability that the earnings growth is above the bound is 90\%.
This suggests that the treatment induces greater reduction in
downside risk but less increase in the upside potential in the
earnings growth. For the right picture the interpretation is just
the opposite.

\begin{figure}[!h]
\centering
 \includegraphics[width=7.8cm, height = 7cm]{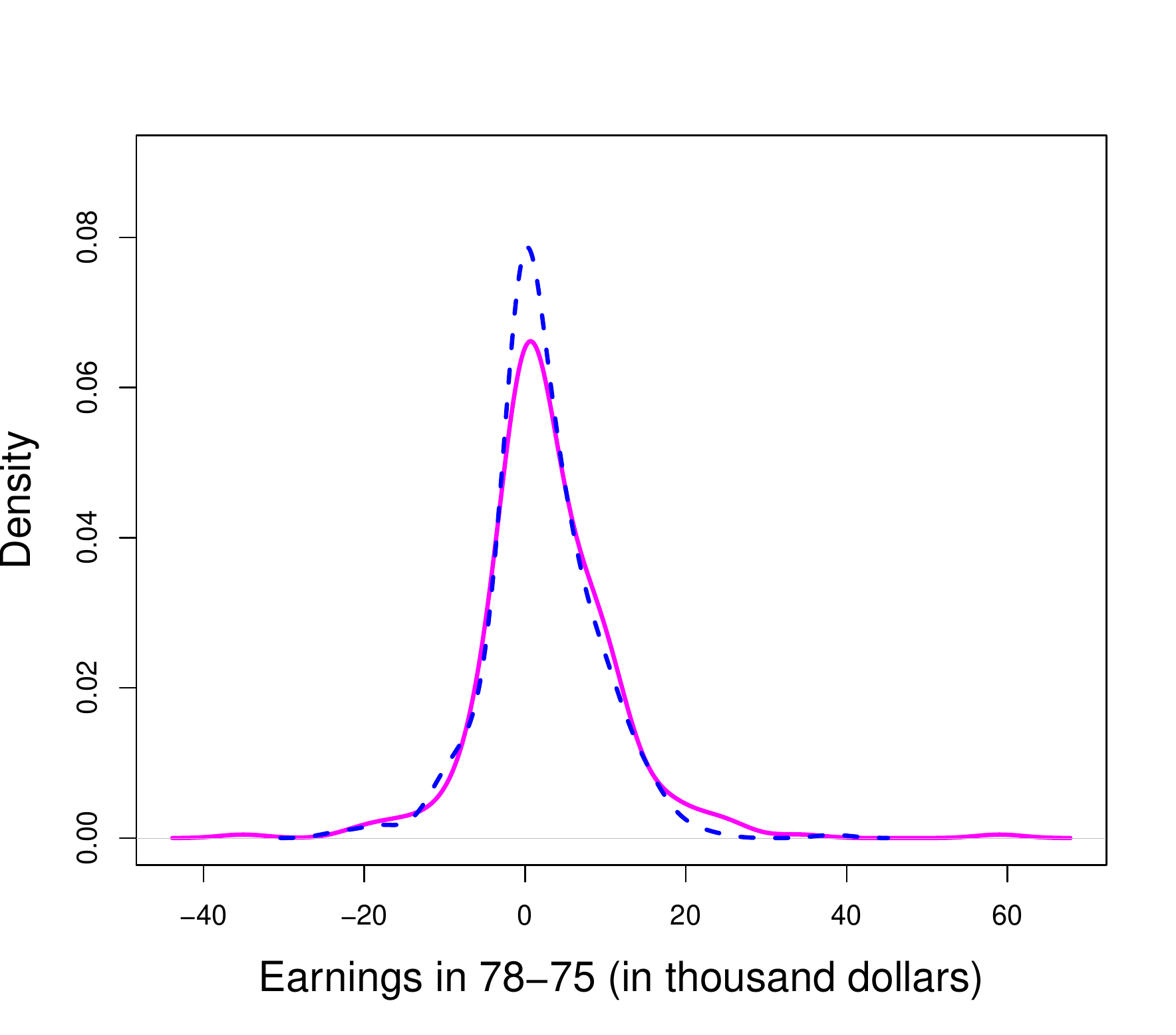}
 \includegraphics[width=7.8cm, height = 7cm]{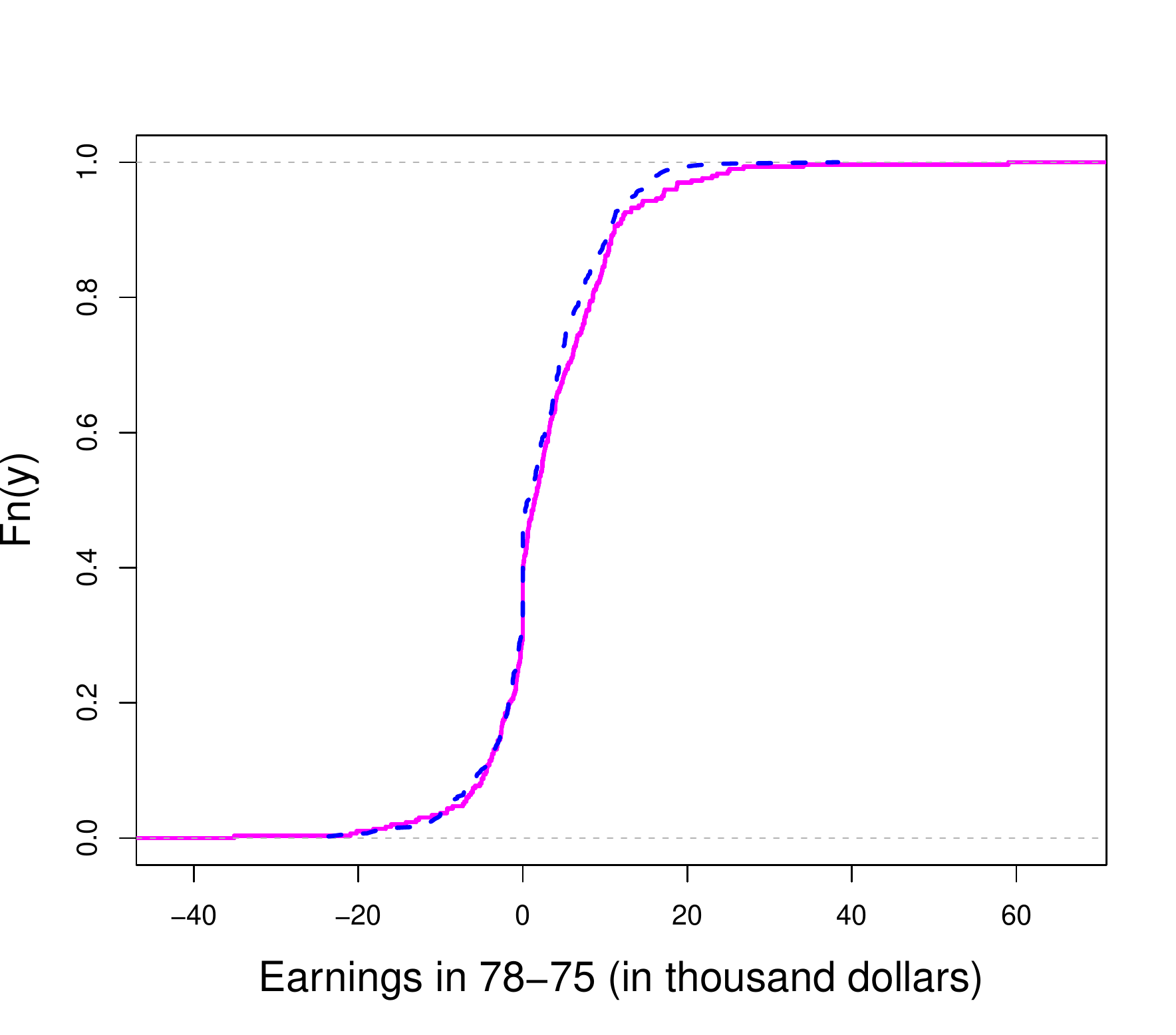}
\caption{Unconditional empirical density function (left) and distribution function (right) of the difference of earnings from 1975 to 1978. The dashed line is associated with the control group and
the solid line is associated with the treatment group.}\label{uncond.dens.y}
\end{figure}

To see which type of stochastic dominance the NSW demonstration
program belongs to, we apply the same data as \cite{DE:2013} for
testing the hypothesis of positive quantile treatment effect for
several quantile levels $\tau$. The data consist of 297 treatment
group observations and 423 control group observations. The response
variable $Y_0$ ($Y_1$) denotes the difference in earnings of control
(treatment) group between 1978 (year of postintervention) and 1975
(year of preintervention). We first apply common statistical
procedures to describe the distribution of these two variables.
Figure \ref{uncond.dens.y} shows the unconditional densities and
distribution function. The cross-validated bandwidth for $\hat
f_0(y)$ is $2.273$ and $2.935$ for $\hat f_1(y)$. The left figure of
Figure \ref{uncond.dens.y} shows the unconditional densities of the
income difference for treatment group and control group. The density
of the treatment group has heavier tails while the density of the
control group is more concentrated around zero. The right
figure shows that 
the two unconditional distribution functions are very close on the
left of the 50\% percentile, and slight deviation appears when the
two distributions are getting closer to 1. Table \ref{quantile.data}
shows that, though the differences are small, but the quantiles of
the unconditional cdf of treatment group are mildly greater than
that of the control group for each chosen $\tau$. The two-sample
Kolmogorov-Smirnov and Cram\'er-von Mises tests, however, yield
results shown in the Table \ref{ecdf.test} which cannot reject the
null hypothesis that the empirical cdfs for
the two groups are the same with confidence levels 1\% or 5\%. 
\begin{table}[h!]
        \begin{center}
        \begin{tabular}{lrrrrrrr}
\hline\hline
$\tau$(\%)           &\multicolumn{1}{c}{10}&\multicolumn{1}{c}{20} &\multicolumn{1}{c}{30} &\multicolumn{1}{c}{50} &\multicolumn{1}{c}{70} &\multicolumn{1}{c}{80} &\multicolumn{1}{c}{90}\\
\hline
Treatment   &-4.38 &-1.55 &0.00 &1.40 &5.48 &8.50 &11.15\\
Control &-4.91 &-1.73 &-0.17 &0.74 &4.44 &7.16 &10.56\\
\hline\hline
        \end{tabular}
        \end{center}
\caption{The unconditional sample quantiles of treatment and control
groups.}\label{quantile.data}
        \end{table}
\begin{table}[h!]
        \begin{center}
        \begin{tabular}{lrr}
\hline\hline
Type of test & Statistics &$p$-value \\
\hline
Kolmogorov-Smirnov &  0.0686 & 0.3835 \\
Cram\'er-von Mises   &  0.2236 & 0.7739 \\
\hline\hline
        \end{tabular}
        \end{center}
\caption{The two sample empirical cdf tests results for treatment
and control groups.}\label{ecdf.test}
        \end{table}

Next we apply our test on quantile regression to evaluate the
treatment effect. In order to compare with \cite{DE:2013}, we first
focus on the case of a one-dimensional covariate. The first
covariate $X_{1i}$ is the age. The second covariate $X_{2i}$ is the
number of years of schooling. The sample values of schooling years
lie in the range of $[3,16]$ and age lies between $[17,55]$. In
order to avoid boundary effect and sparsity of the samples, we look
at the ranges [7,13] for schooling years and [19,31] for age. We
apply the bootstrap CC method for quantiles
$\tau=0.1,0.2,0.3,0.5,0.7,0.8$ and $0.9$. We apply the quartic
kernel. The cross-validated bandwidths are chosen in the same way as
for conditional densities with the \texttt{R} package \texttt{np}.
The resulting bandwidths are (2.2691,2.5016) for the treatment group
and (2.7204, 5.9408) for the control group. In particular, for
smoothing the data of the treatment group, for $\tau=0.1$ and $0.9$,
we enlarge the cross-validated bandwidths by a constant of $1.7$;
for $\tau=0.2,0.3,0.7,0.8$, the cross-validated bandwidths are
enlarged by constant factor $1.3$. These inflated bandwidths are
used to handle violent roughness in extreme quantile levels. The
bootstrap CCs are computed with 10,000 repetitions. The level of the
test is $\alpha=5\%$.

\begin{figure}[!h]
\centering
\includegraphics[width=5cm, height = 5cm]{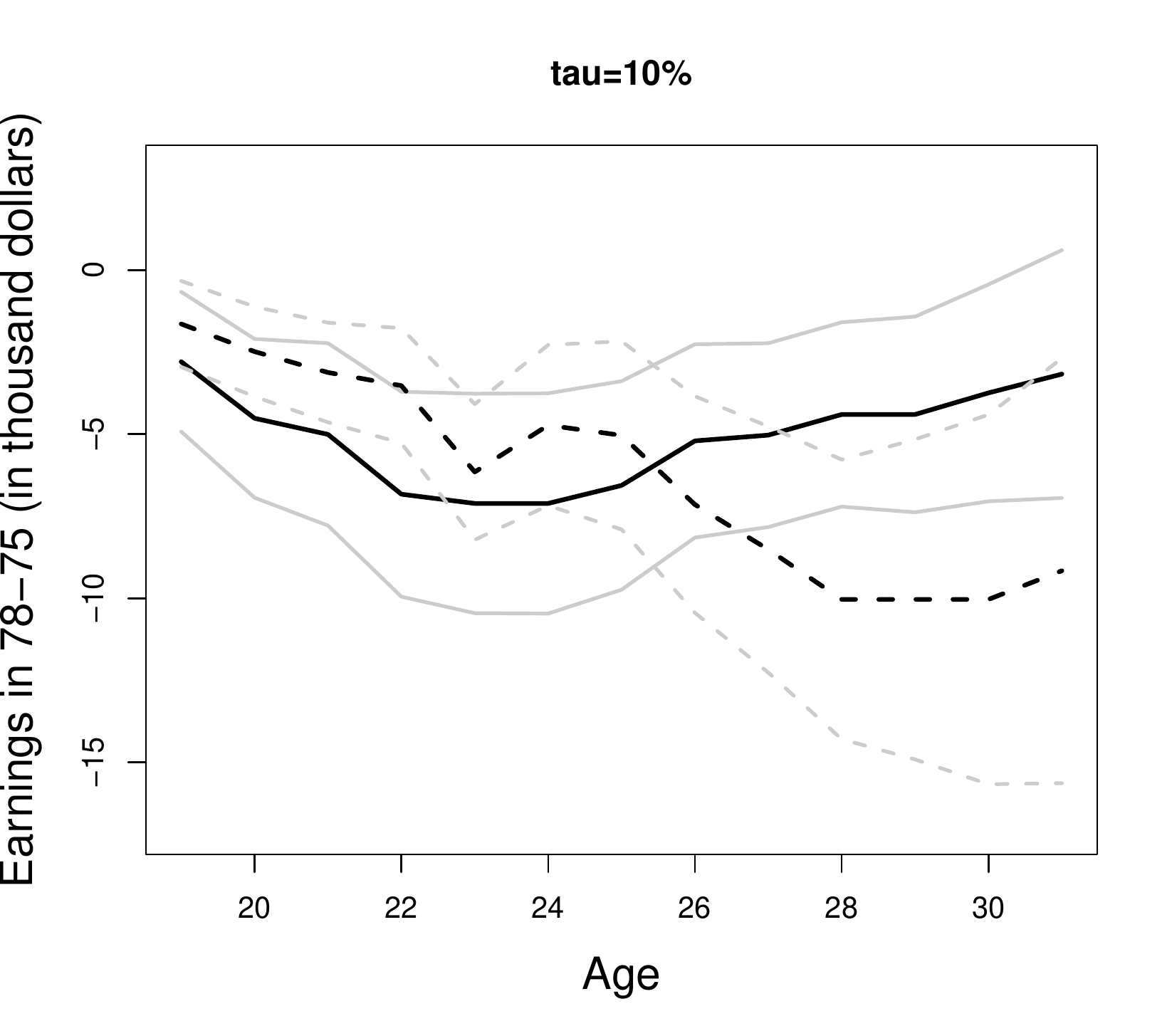}
\includegraphics[width=5cm, height = 5cm]{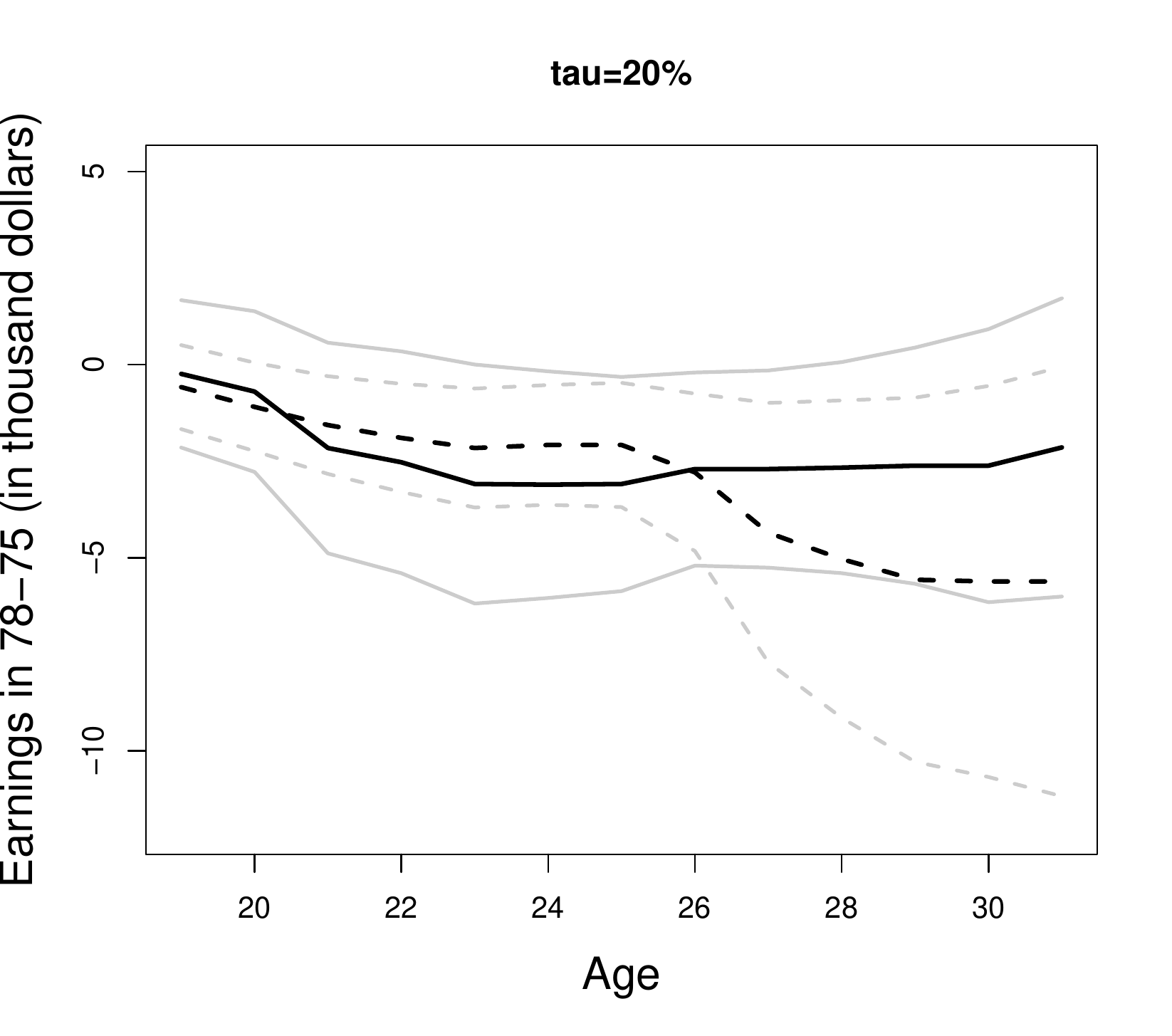}
\includegraphics[width=5cm, height = 5cm]{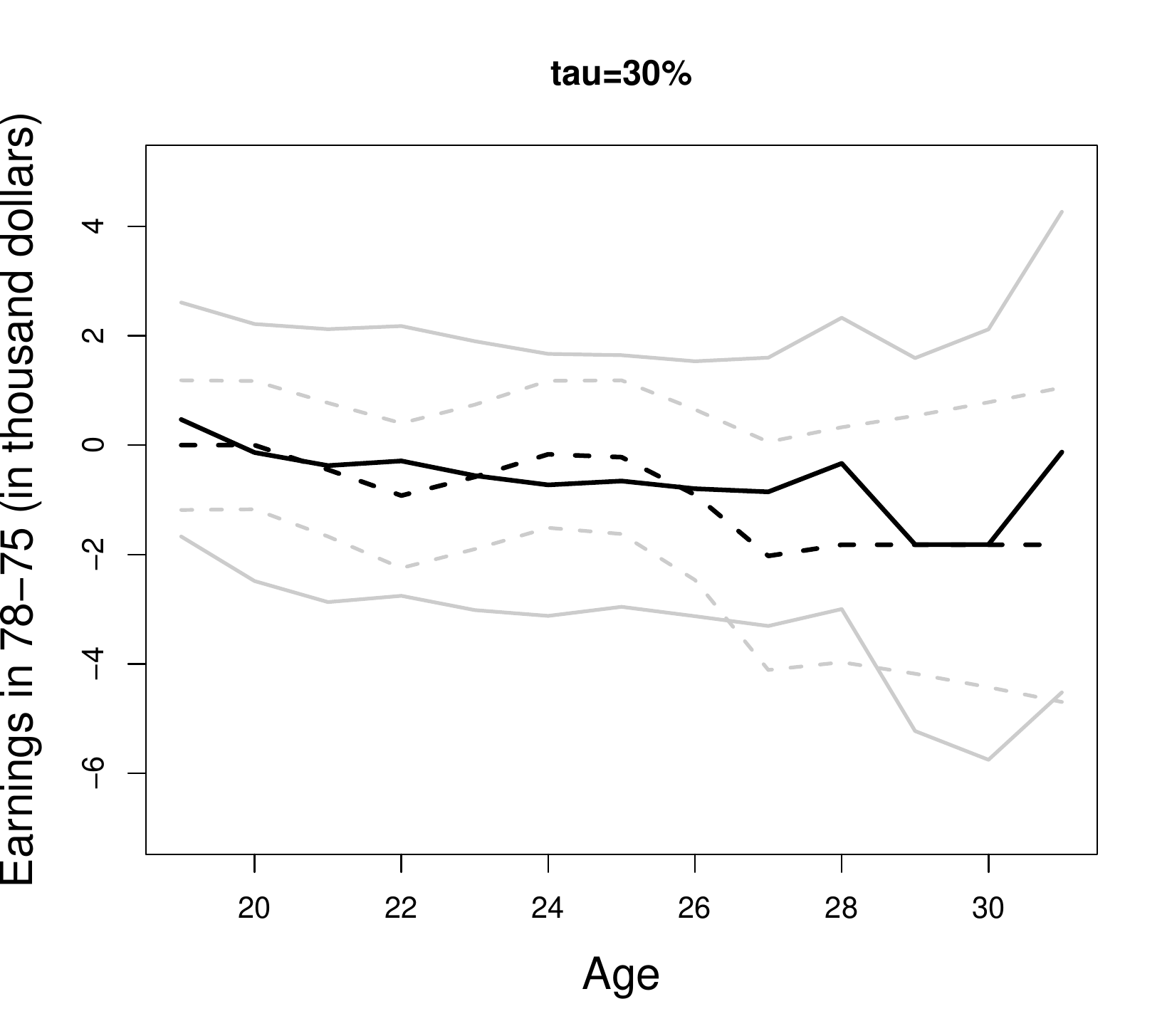}
\includegraphics[width=5cm, height = 5cm]{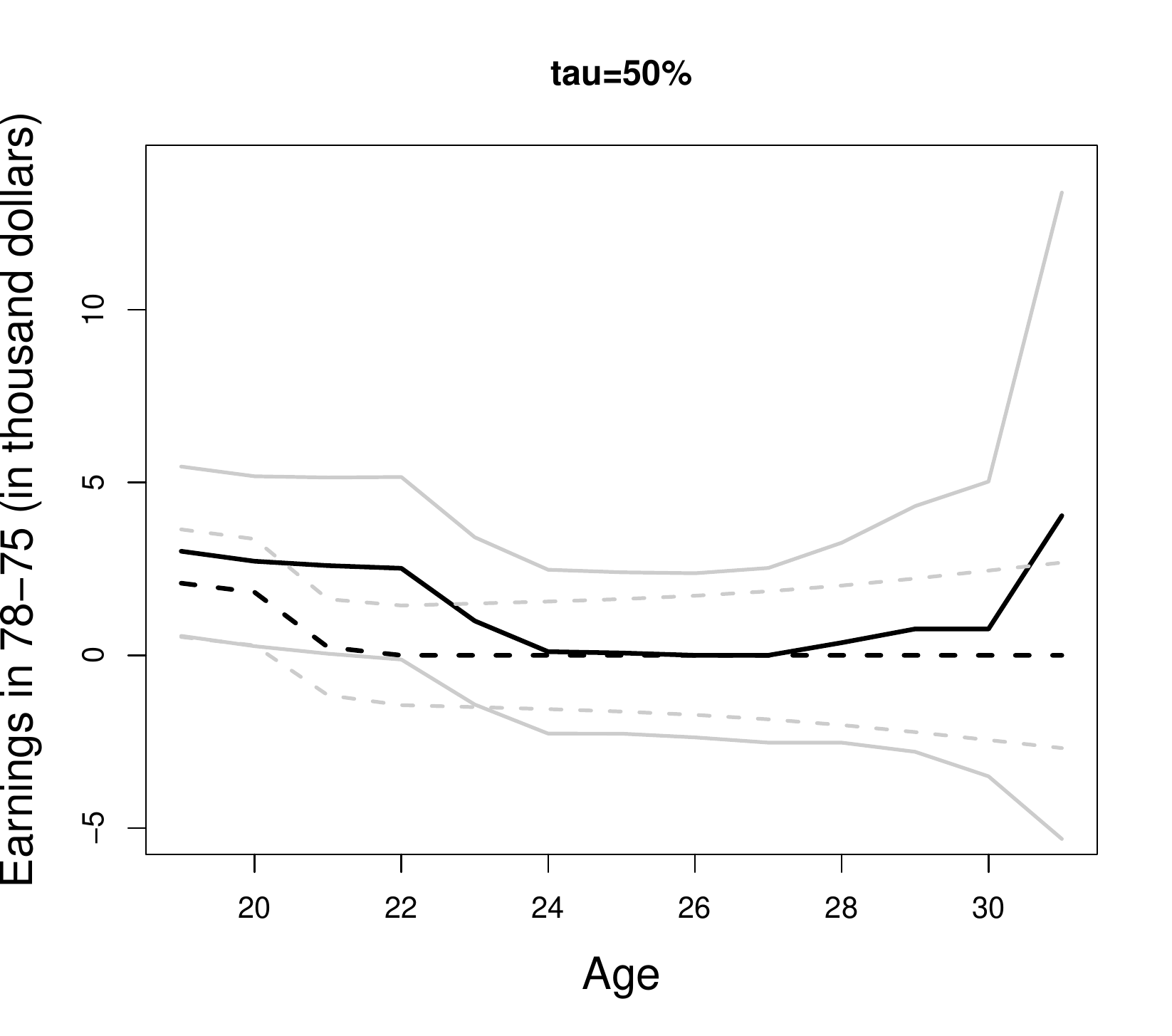}
\includegraphics[width=5cm, height = 5cm]{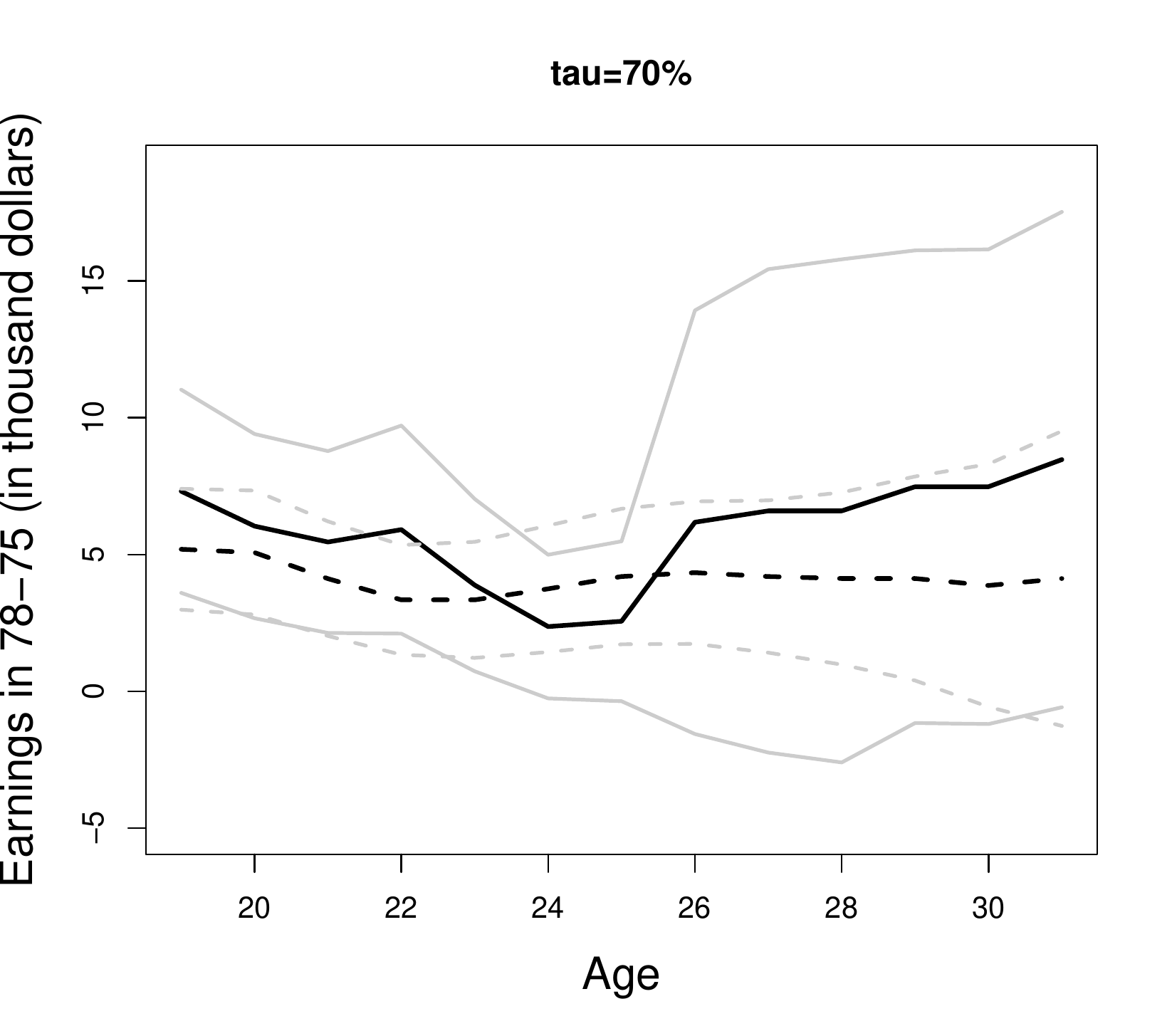}
\includegraphics[width=5cm, height = 5cm]{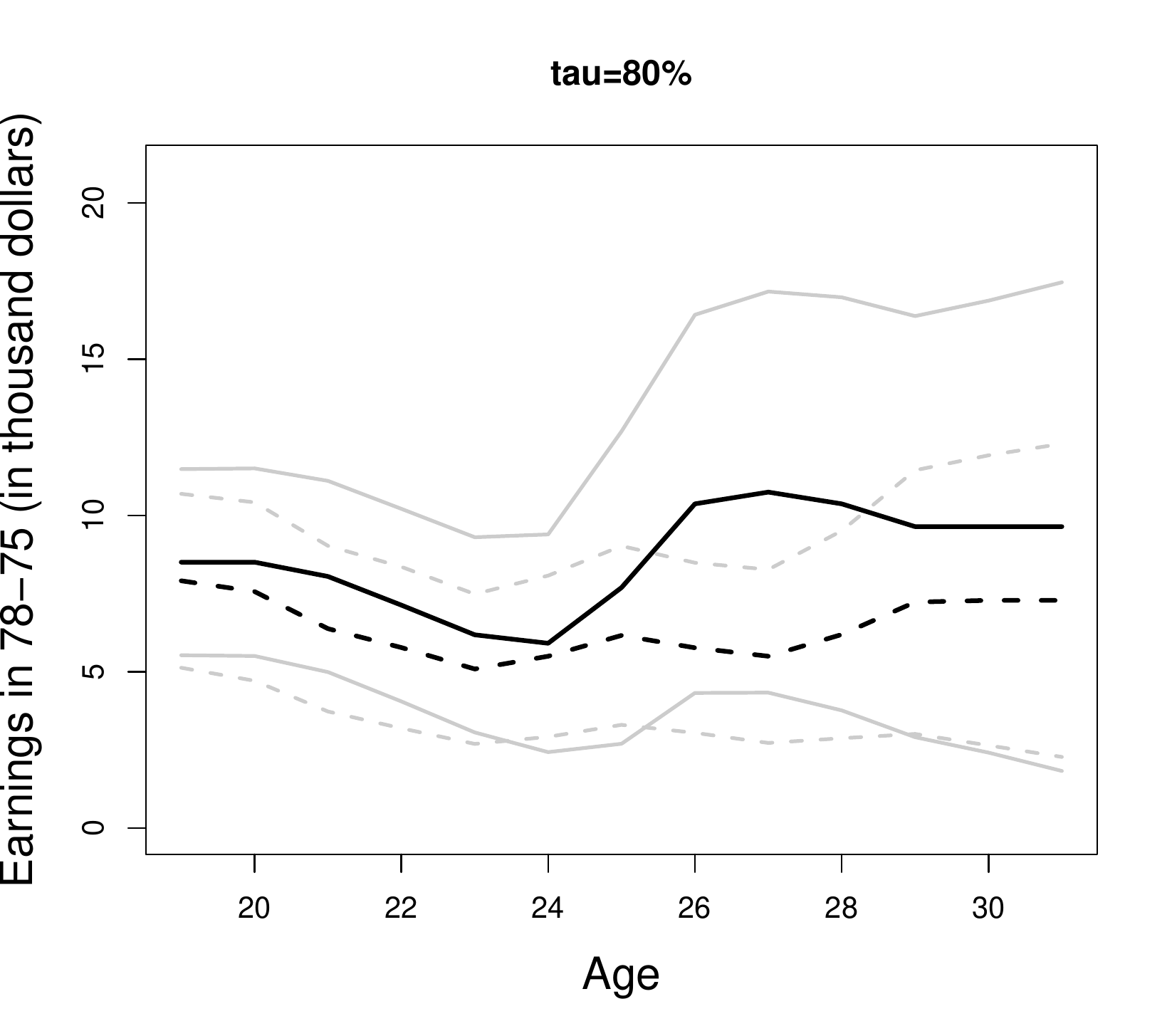}
\includegraphics[width=5cm, height = 5cm]{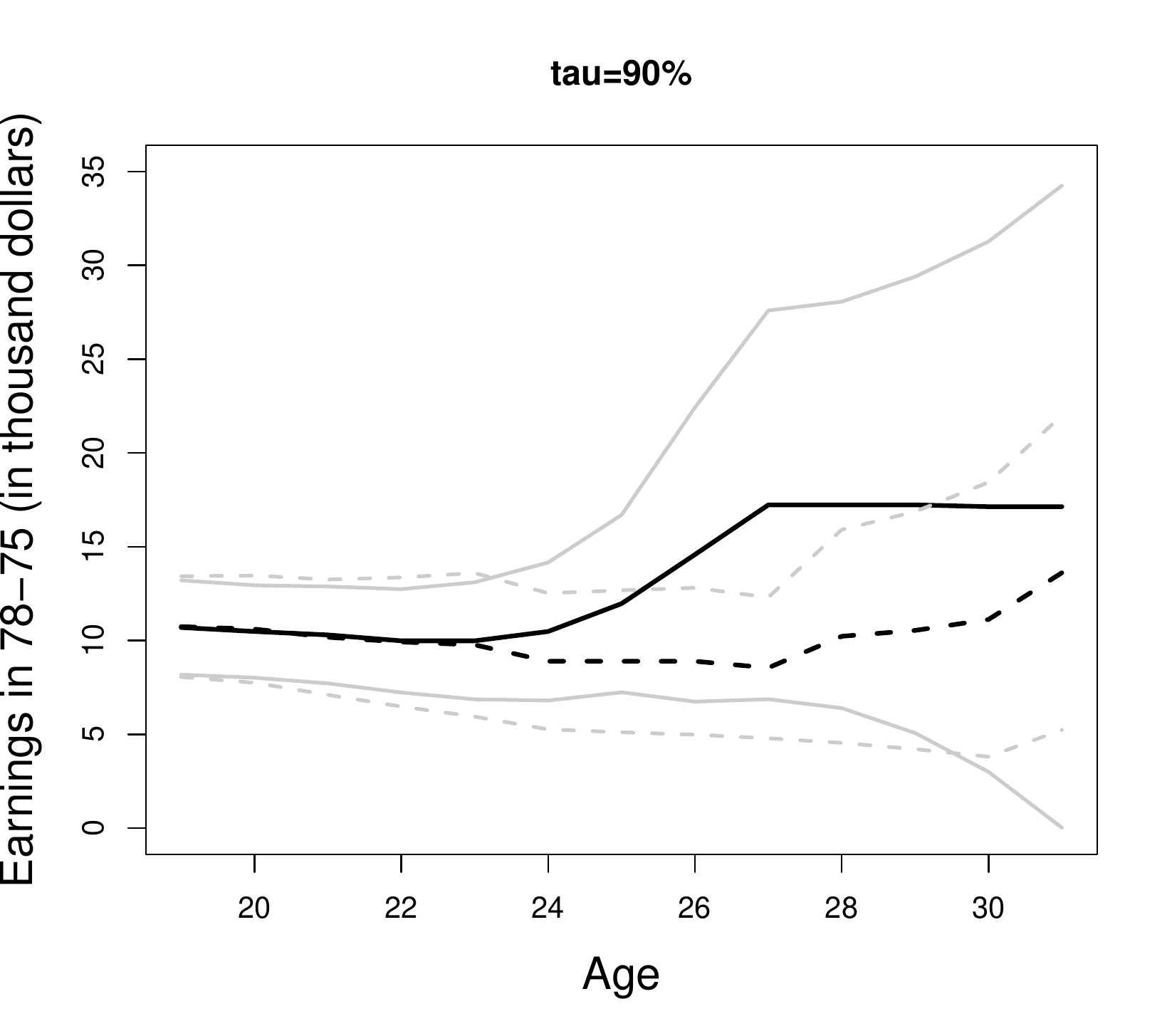}
\caption{Nonparametric quantile regression estimates and CCs for the changes in earnings between 1975-1978 as a function of age. The solid dark lines correspond to the conditional quantile of the treatment group and the solid light lines sandwich its CC, and the dashed
dark lines correspond to the conditional quantiles of the control group and the solid light lines sandwich its CC.}\label{lalonde.age}
\end{figure}

\begin{figure}[!h]
\centering
\includegraphics[width=5cm, height = 5cm]{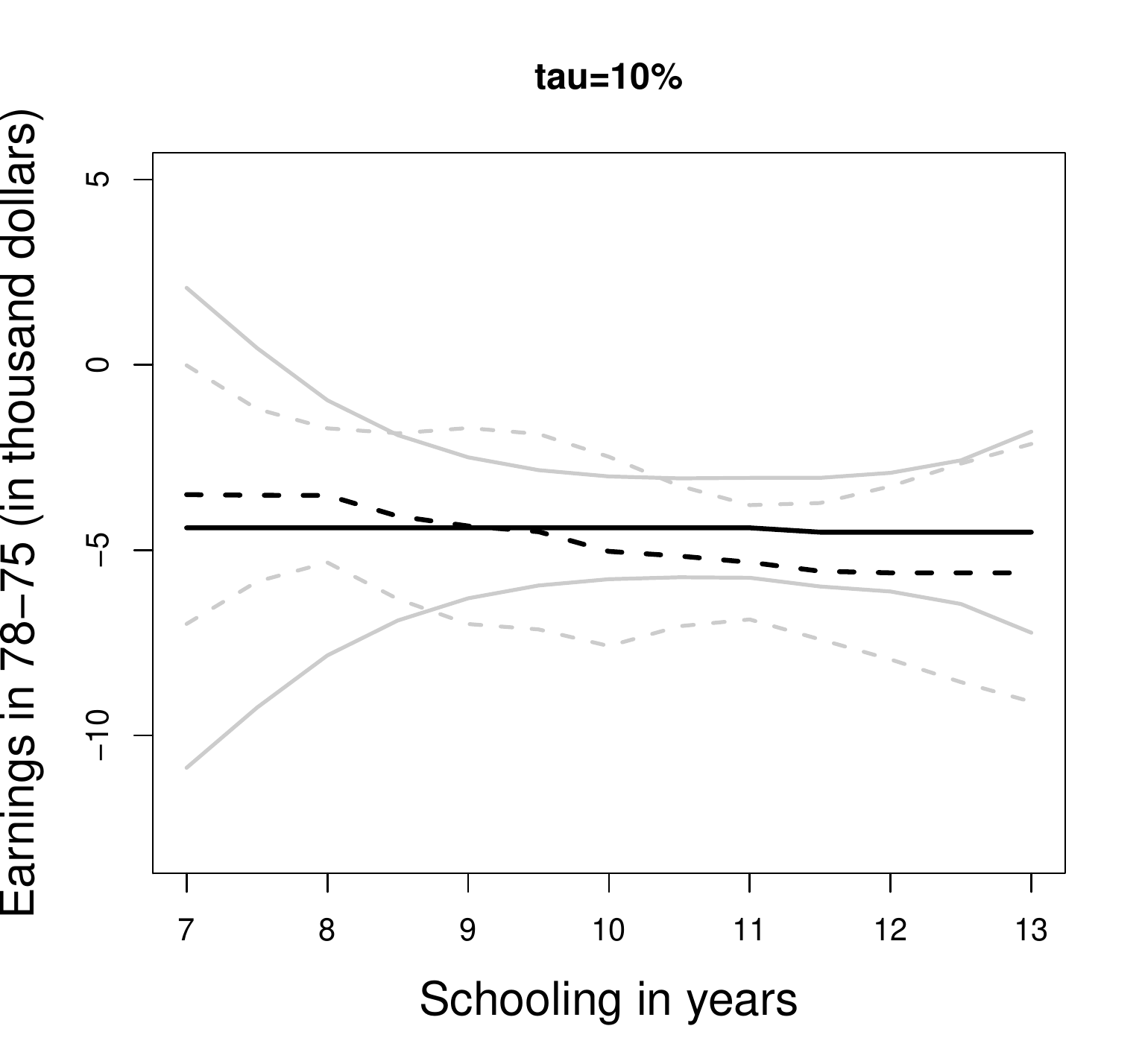}
  \includegraphics[width=5cm, height = 5cm]{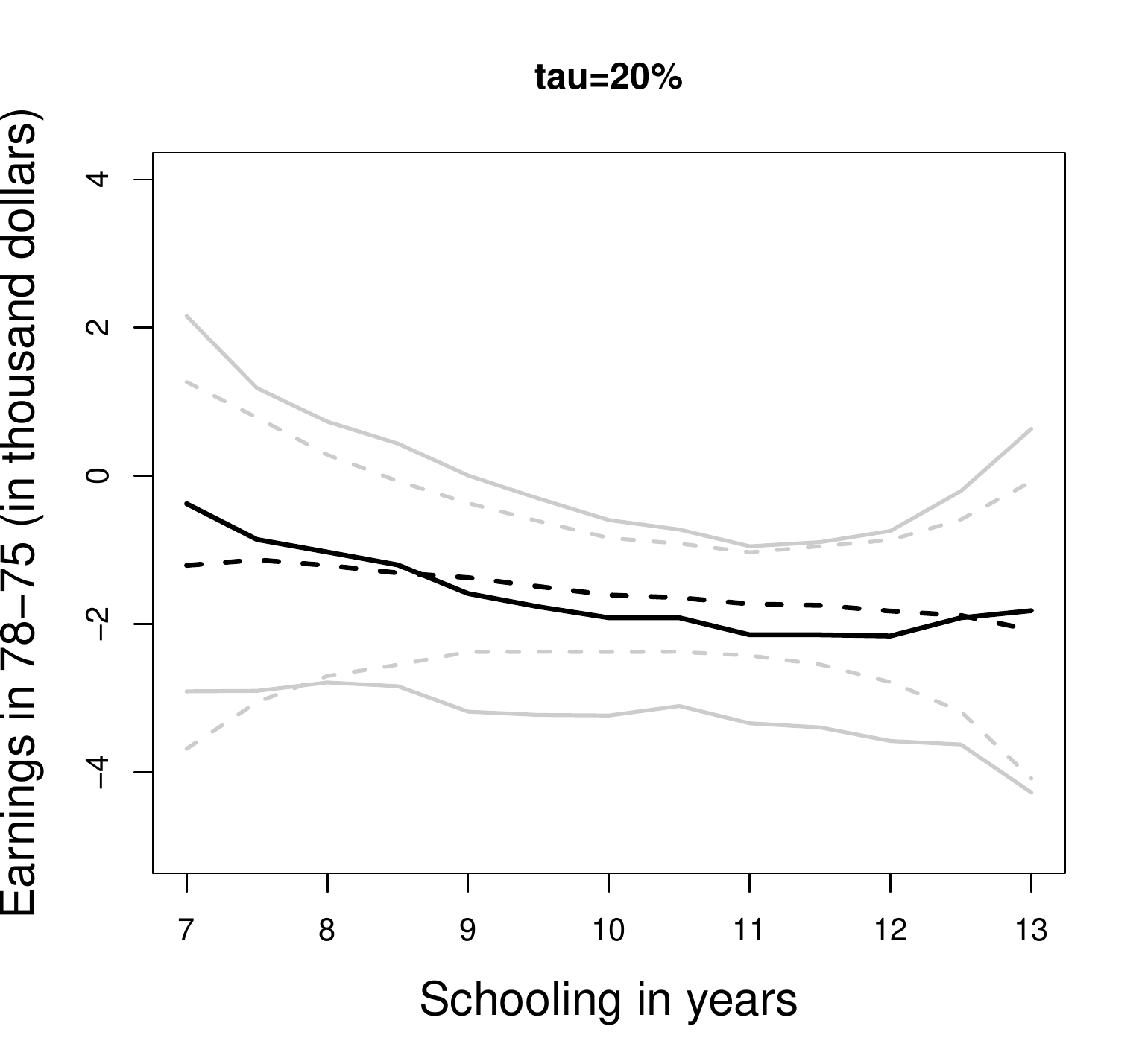}
  \includegraphics[width=5cm, height = 5cm]{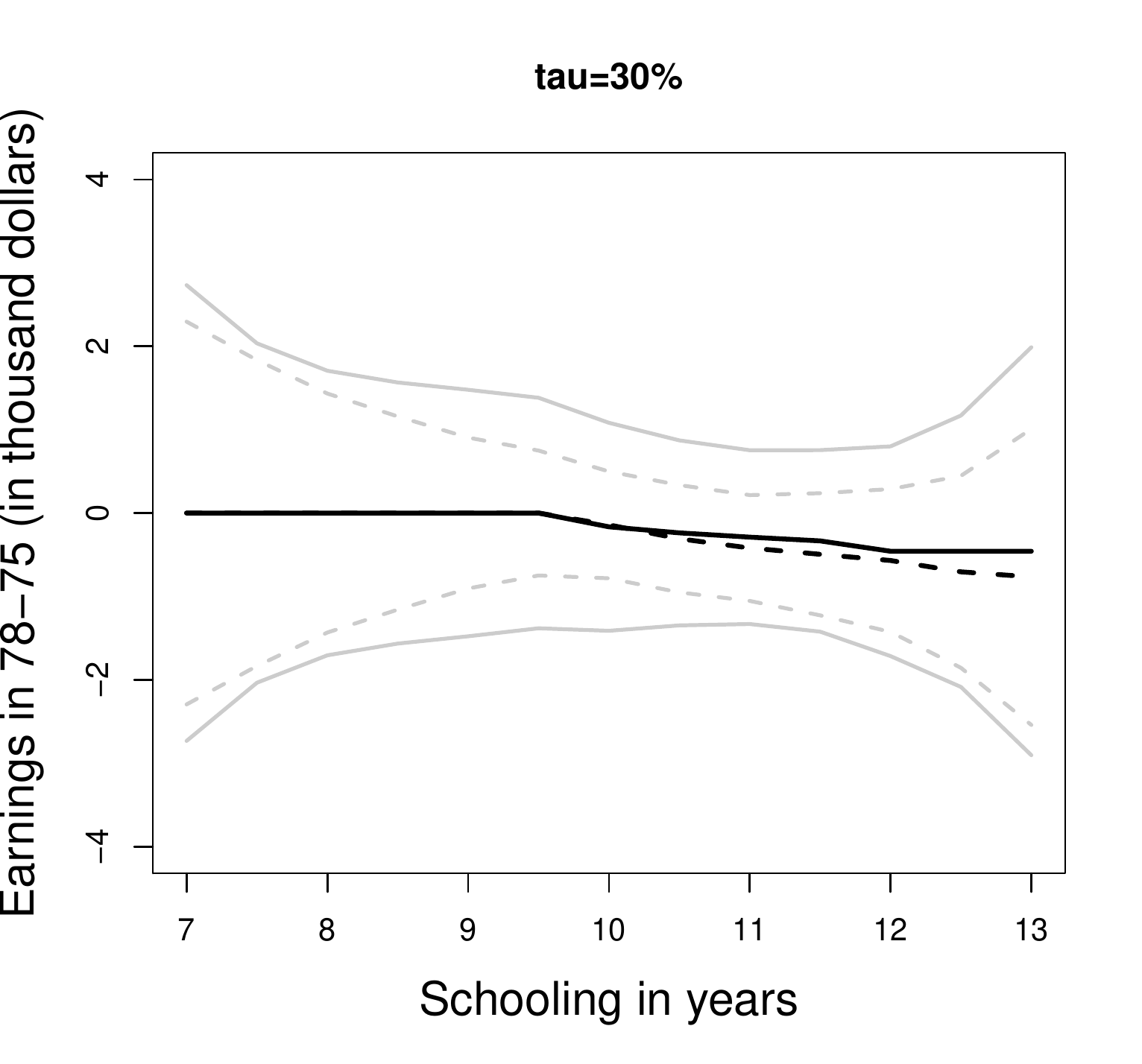}
 \includegraphics[width=5cm, height = 5cm]{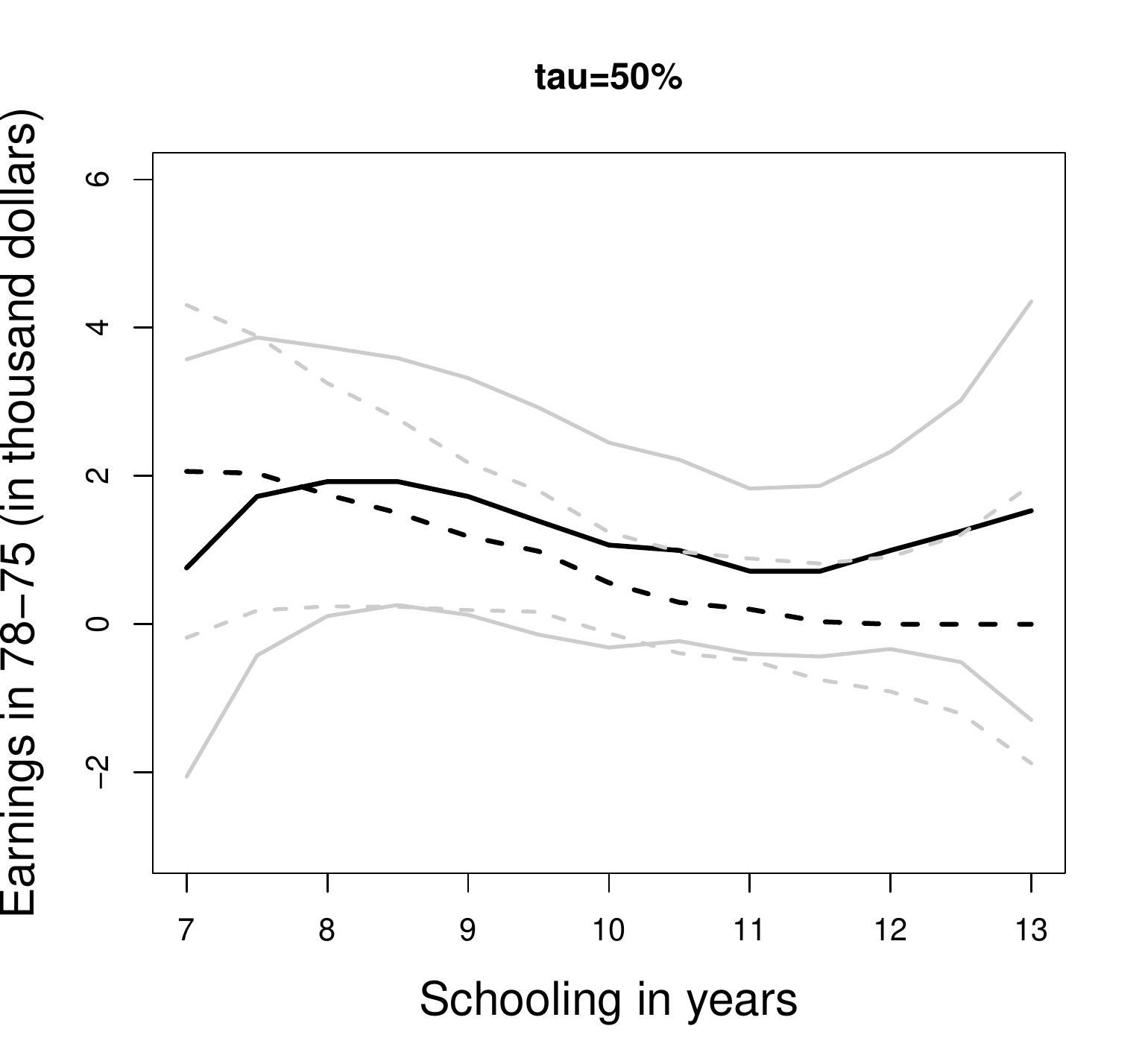}
  \includegraphics[width=5cm, height = 5cm]{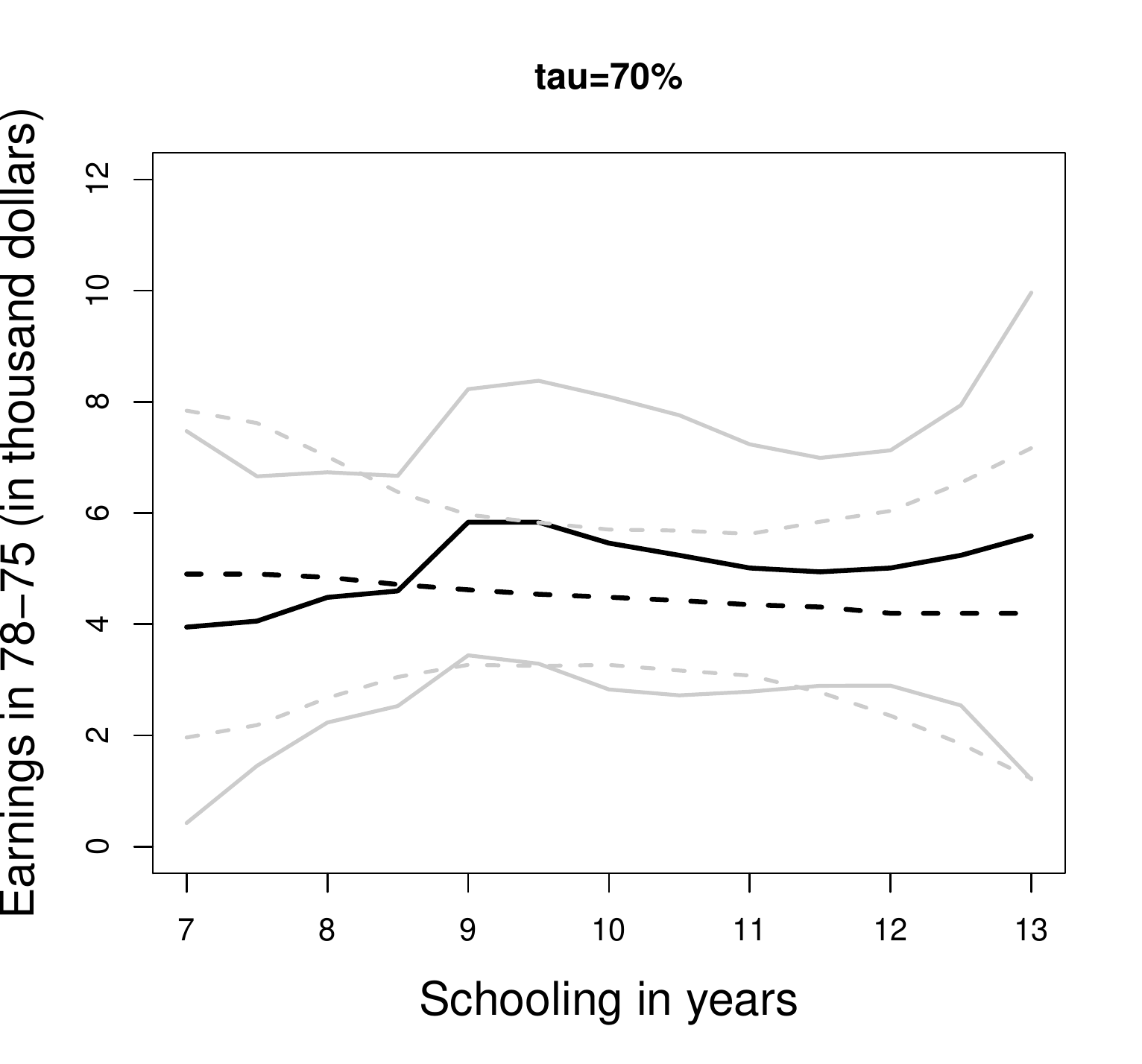}
  \includegraphics[width=5cm, height = 5cm]{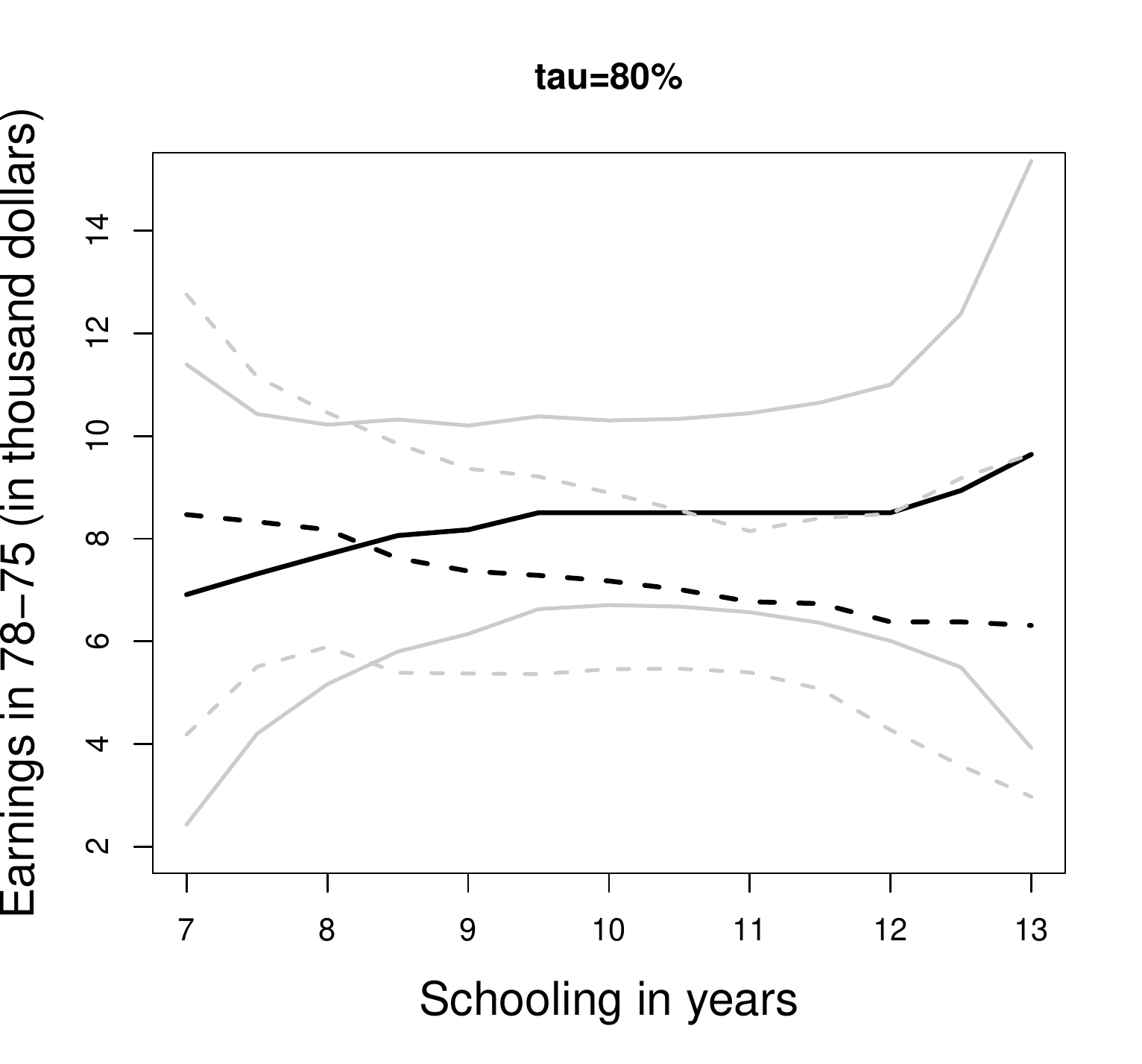}
  \includegraphics[width=5cm, height = 5cm]{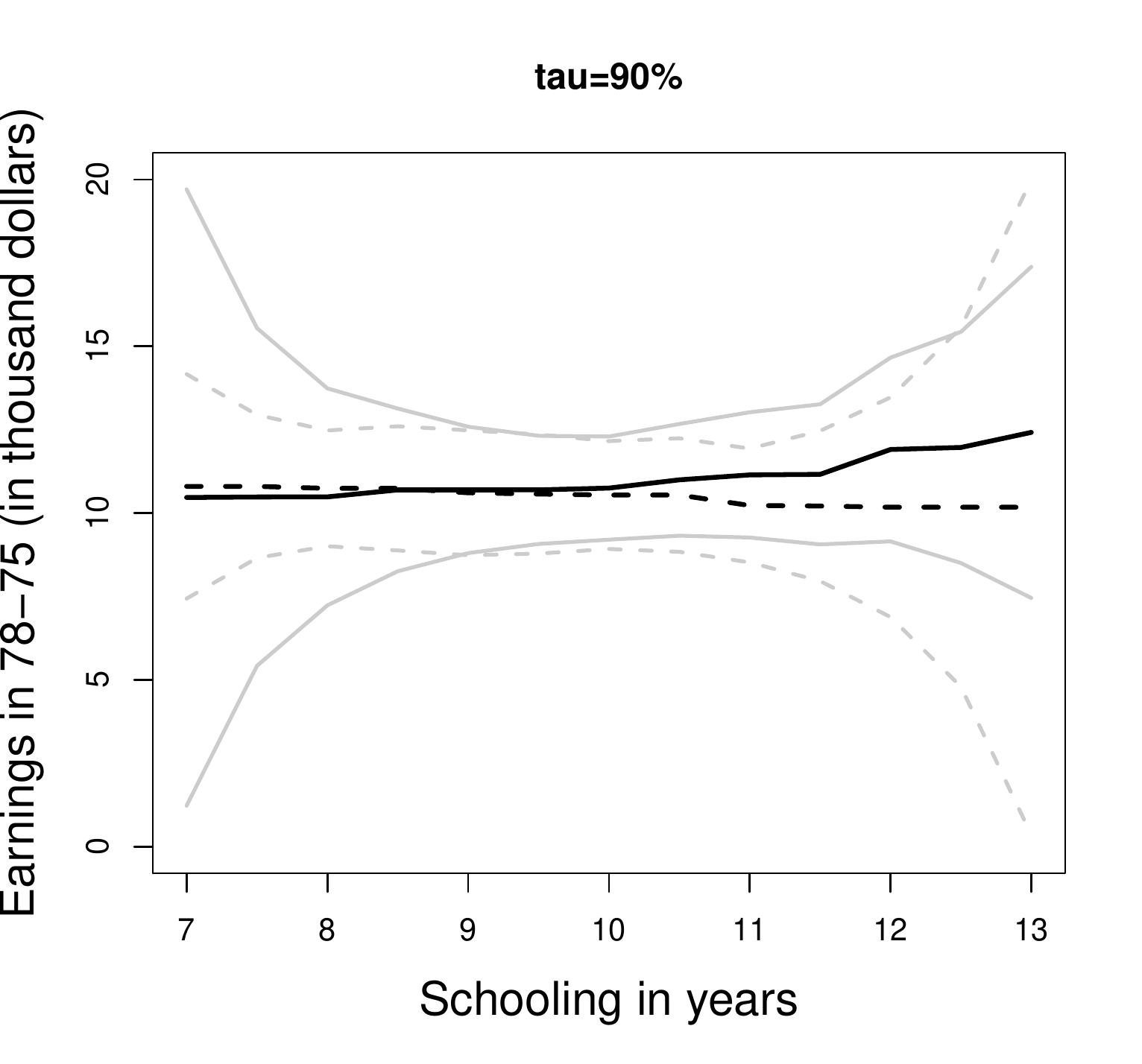}
\caption{Nonparametric quantile regression estimates and CCs for the changes in earnings between 1975-1978 as a function of years of schooling. The solid dark lines correspond to the conditional quantile of the treatment group and the solid light lines sandwich its CC, and the dashed
dark lines correspond to the conditional quantiles of the control group and the solid light lines sandwich its CC.}\label{lalonde.educ}
\end{figure}

The results of the two quantile regressions with one-dimensional
covariate, and their CCs for various quantile levels are presented
in Figure \ref{lalonde.age} and \ref{lalonde.educ}. We observe that
for all chosen quantile levels the quantile estimates associated to
the treatment group lie above that of the control group when age is
over certain levels, and particularly for $\tau=10\%, 50\%, 80\%$
and $90\%$, the quantile estimates for treatment group exceeds the
upper CCs for the quantile estimates of the control group. On the
other hand, at $\tau=10\%$, the quantile estimates for the control
group drop below the CC for treatment group for age greater than 27.
Hence, the results here show a tendency that both the downside risk
reduction and the upside potential enhancement of earnings growth
are achieved, as the older individuals benefit the most from the
treatment. Note that we observe a heterogeneous treatment effect in
age and the weak dominance of the conditional quantiles of the
treatment group with respect to those of the control group, i.e.,
\eqref{csd_quantile} holds for the chosen quantile levels, which are
in line with the findings of \cite{DE:2013}.

We now turn to Figure \ref{lalonde.educ}, where the covariate is the
years of schooling. The treatment effect is not significant for
conditional quantiles at levels $\tau=10\%, 20\%$ and $30\%$. This
suggests that the treatment does little to reduce the downside risk
of the earnings growth for individuals with various degrees of
education. Nonetheless, we constantly observe that the regression
curves of the treatment group rise above that of the control group
after a certain level of the years of schooling for quantile levels
$\tau=50\%, 70\%, 80\%$ and $90\%$. Notice that for $\tau=50\%$ and
$80\%$ the regression curves associated to the treatment group reach
the upper boundary of the CC of the control group. This suggests
that the treatment effect tends to raise the upside potential of the
earnings growth, in particular for those individuals who spent more
years in the school. It is worth noting that we also see a
heterogeneous treatment effect in schooling years, although the
heterogeneity in education is less strong than the heterogeneity in
age.

\begin{figure}[!h]
\centering
\subfigure[$\tau=10\%$]{\includegraphics[width=5cm, height = 5cm]{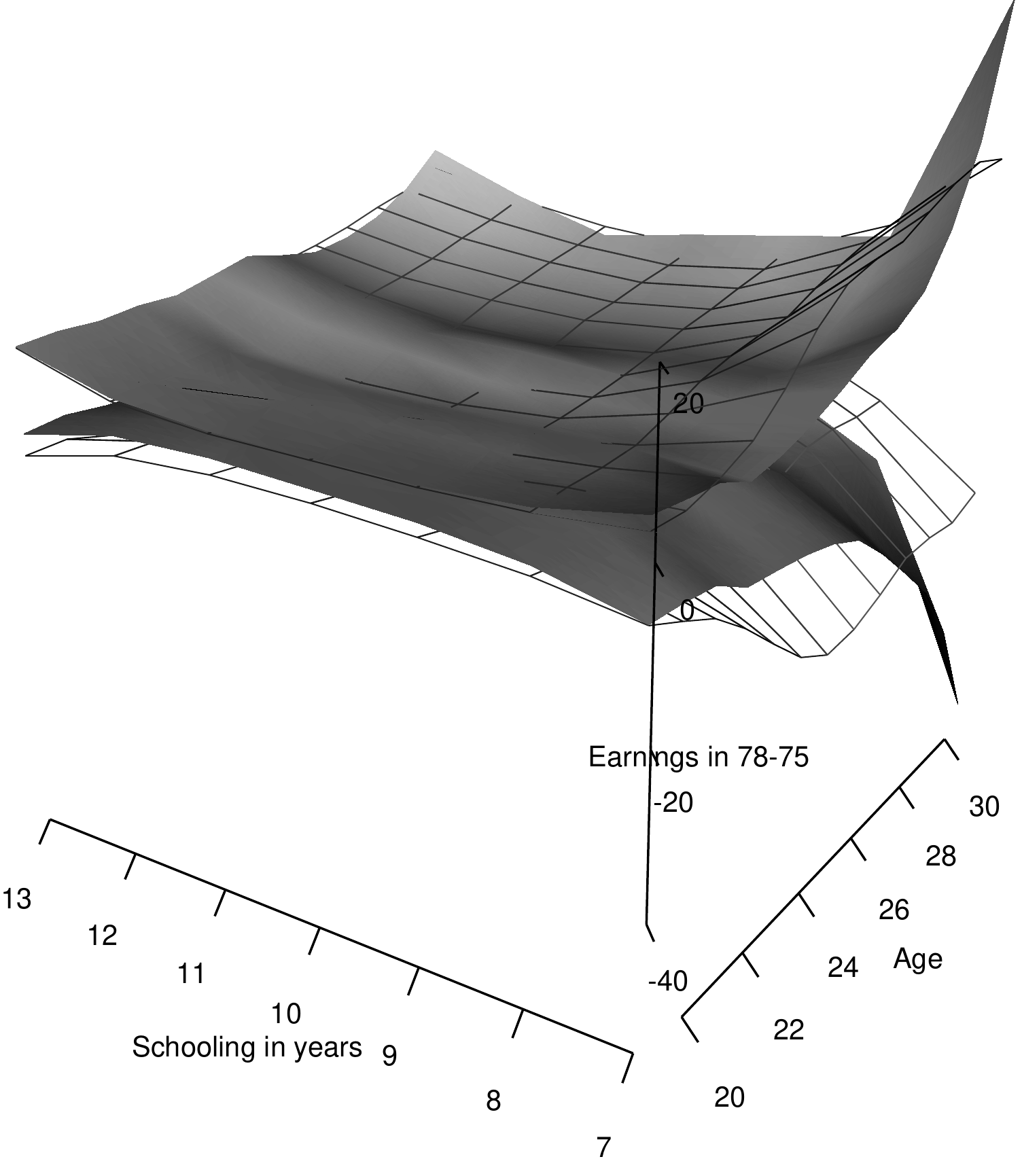}}
\subfigure[$\tau=20\%$]{\includegraphics[width=5cm, height = 5cm]{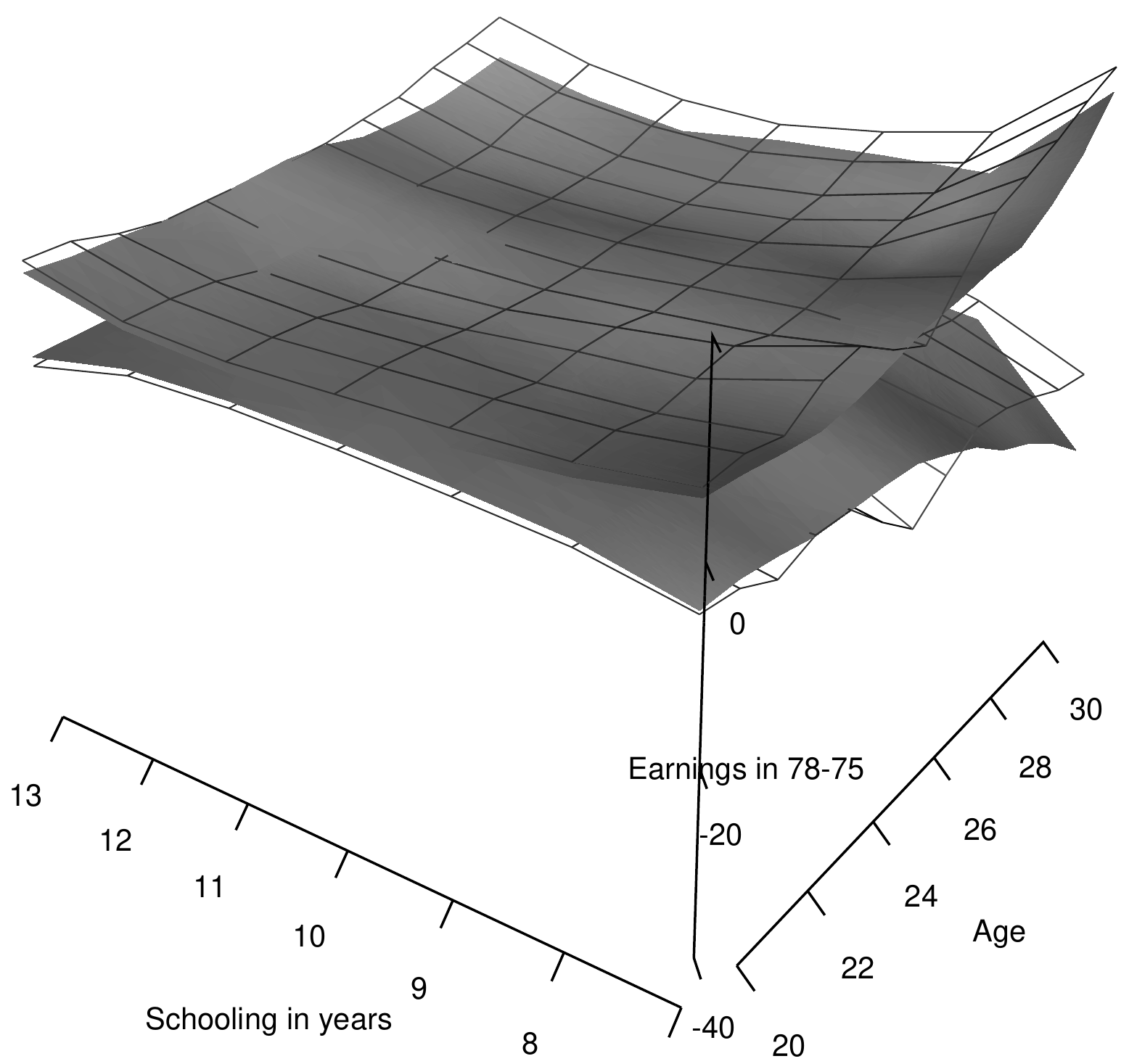}}
\subfigure[$\tau=30\%$]{\includegraphics[width=5cm, height = 5cm]{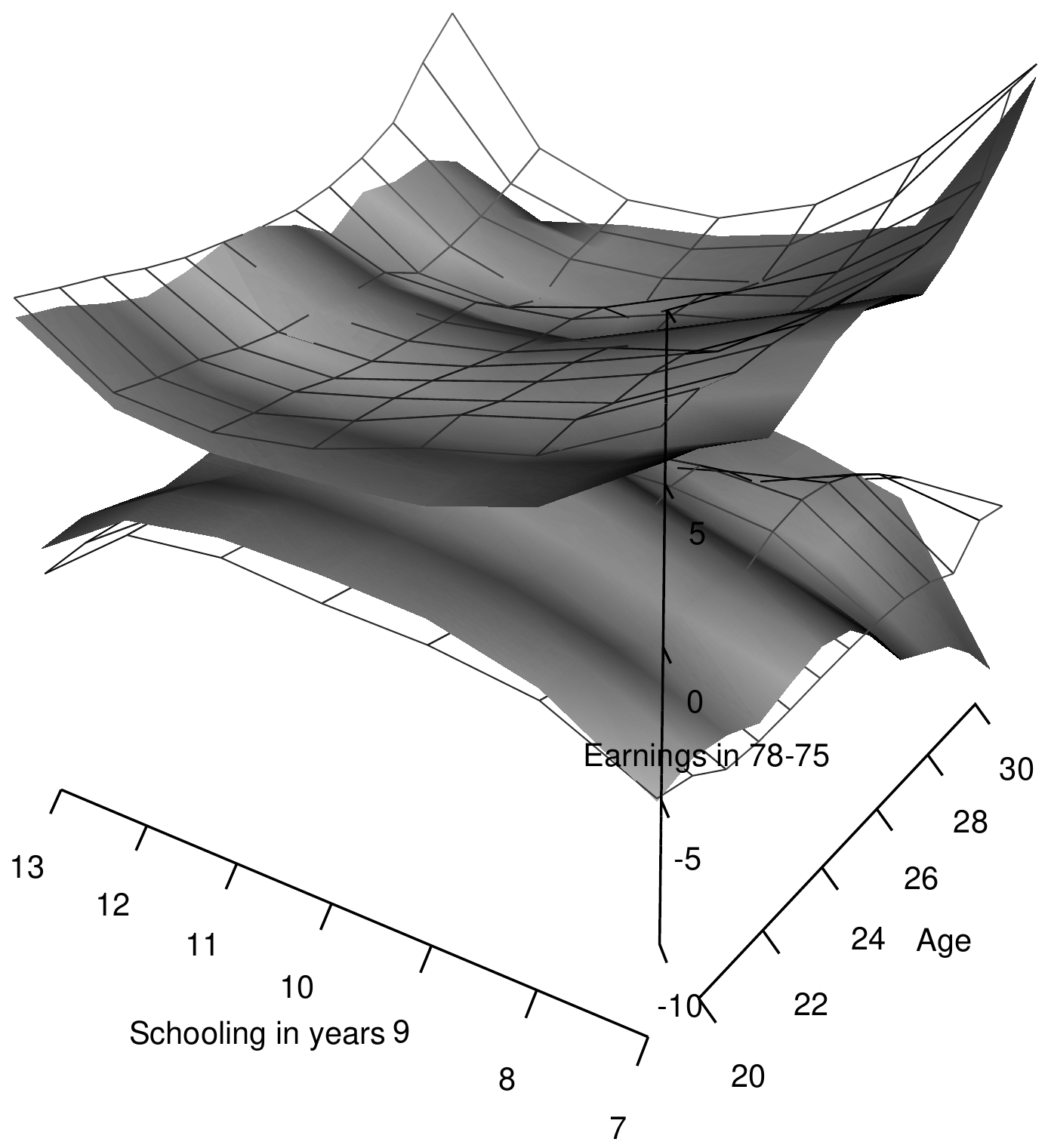}}
\subfigure[$\tau=50\%$]{\includegraphics[width=5cm, height = 5cm]{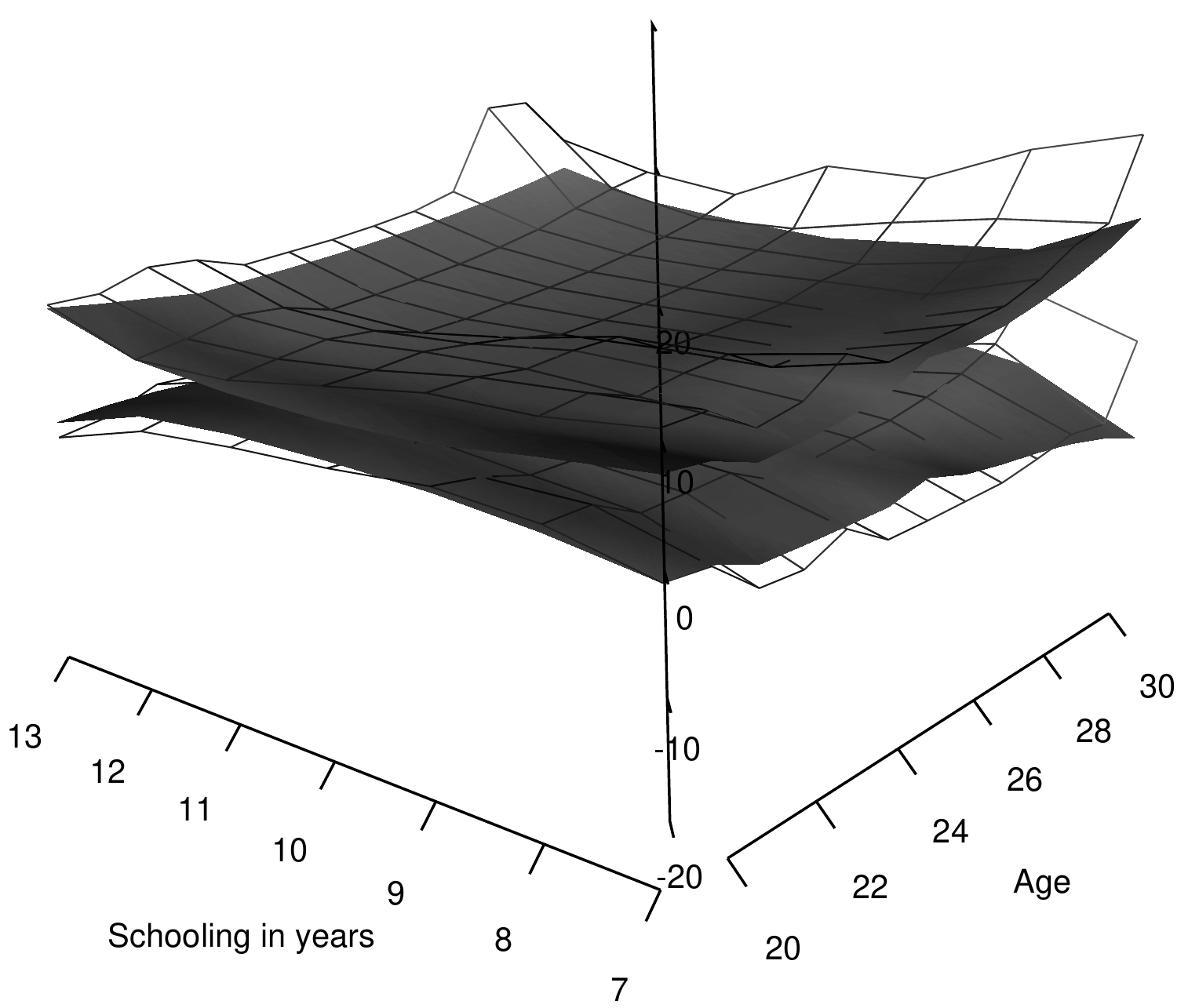}}
\subfigure[$\tau=70\%$]{\includegraphics[width=5cm, height = 5cm]{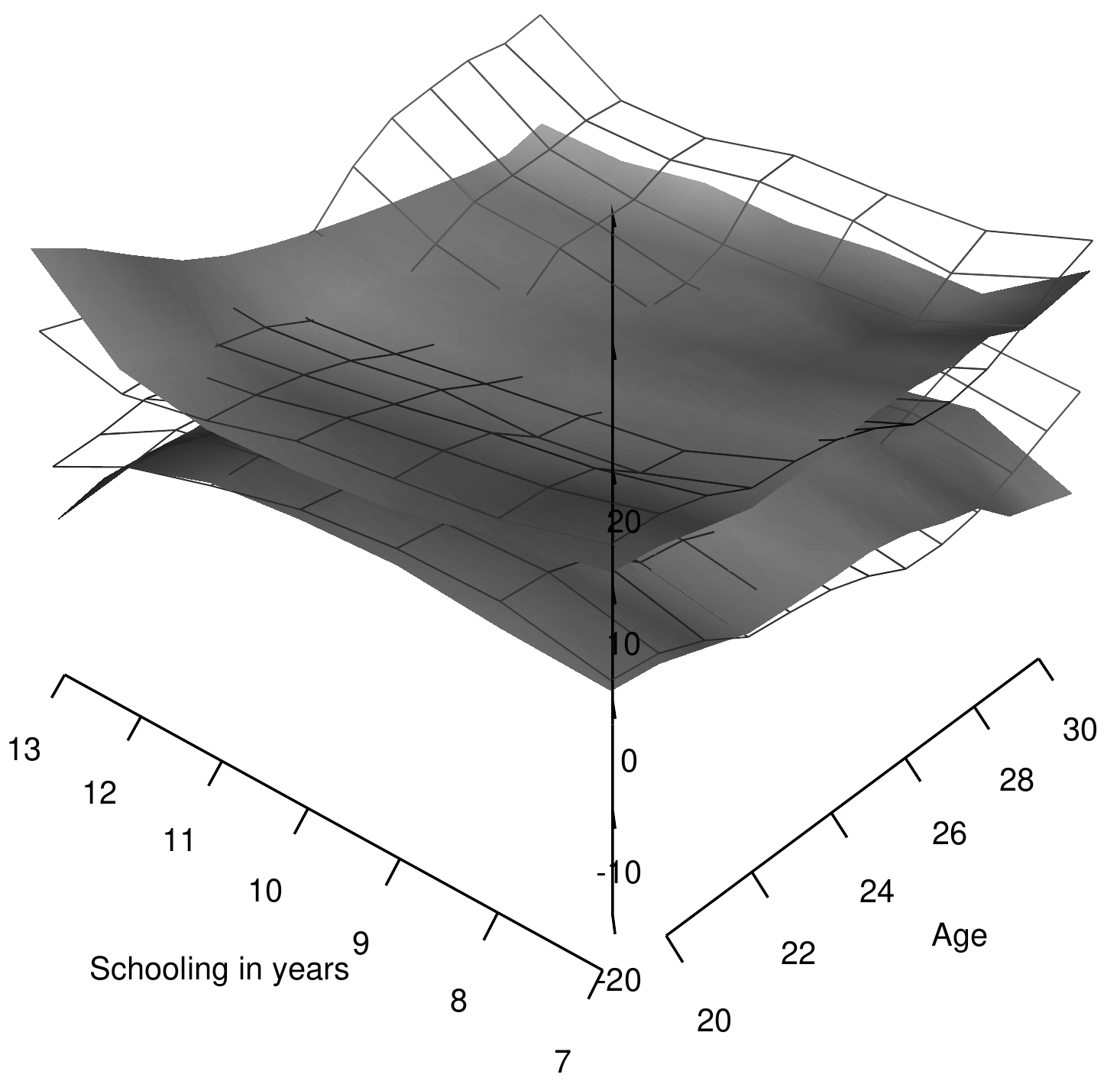}}
\subfigure[$\tau=80\%$]{\includegraphics[width=5cm, height = 5cm]{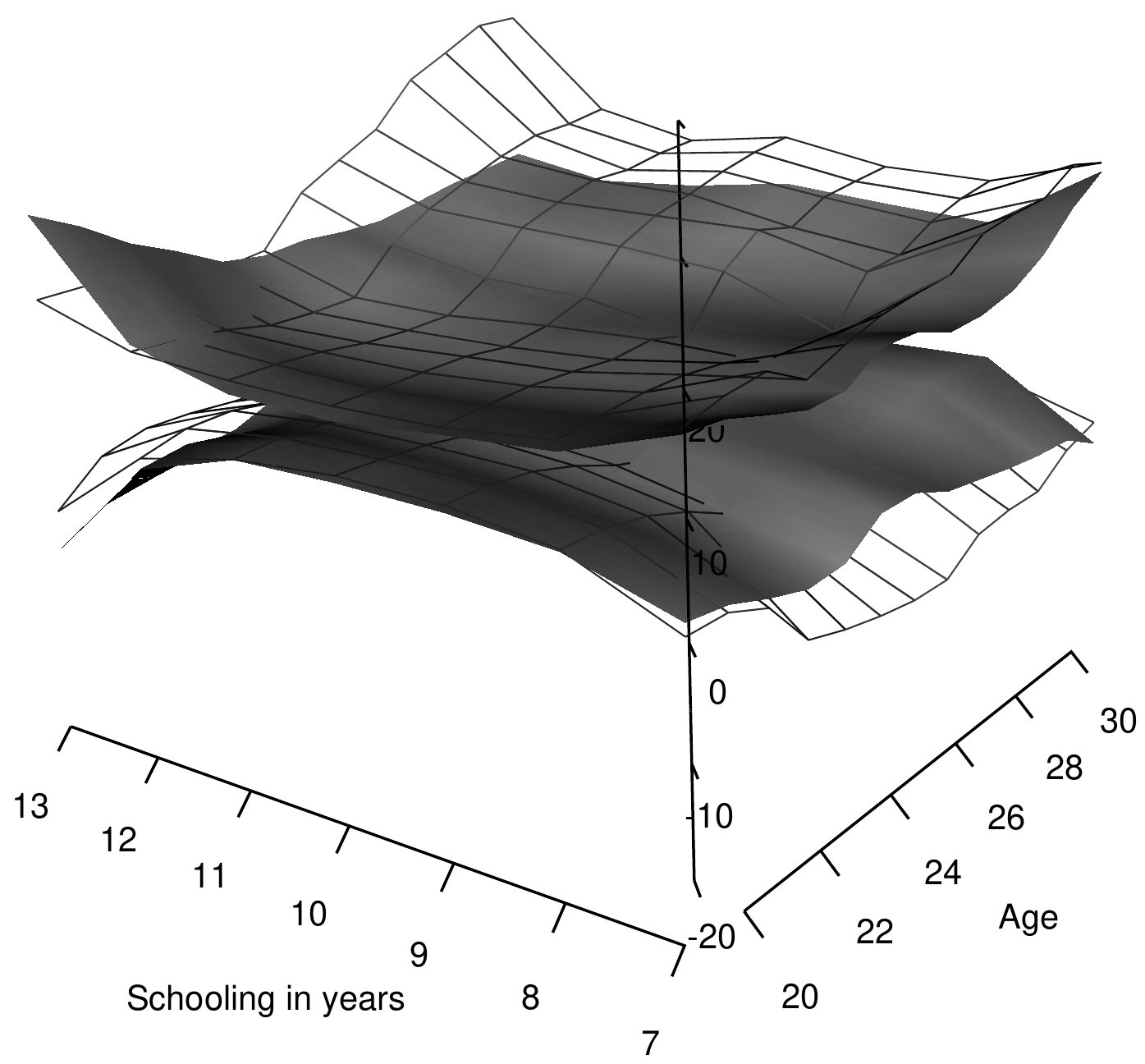}}
\subfigure[$\tau=90\%$]{\includegraphics[width=5cm, height = 5cm]{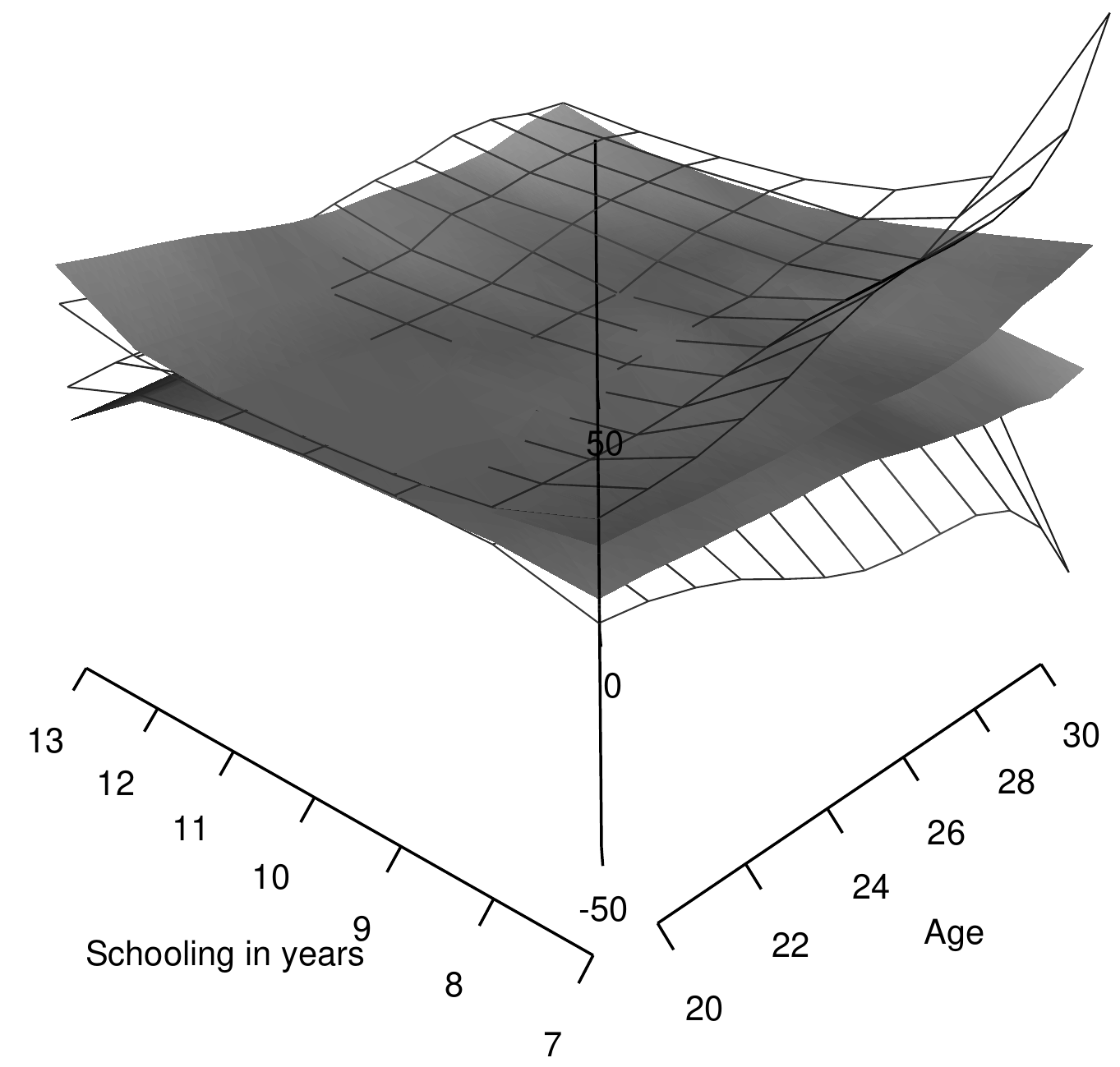}}
\caption{The CCs for the treatment group and the control group. The net surface corresponds to the control group quantile CC and the solid surface corresponds to the treatment group quantile CC.}\label{lalonde.boot}
\end{figure}

The previous regression analyses separately conditioning on
covariates age and schooling years only give a limited view on the
performance of the program, we now proceed to the analysis
conditioning on the
two covariates jointly $(X_{1i},X_{2i})$. The estimation settings are similar to the case of univariate covariate. 
Figure \ref{lalonde.boot} shows the quantile regression CCs. From a
first glance of the pictures, the $\tau$-quantile CC of the
treatment group and that of the control group overlap extensively
for all $\tau$. We could not find sufficient evidence to reject the
null hypothesis that the conditional distribution of treatment group
and control group are equivalent.

\begin{figure}[!h]
\centering
\subfigure[$\tau=10\%$]{\includegraphics[width=5cm, height = 5cm]{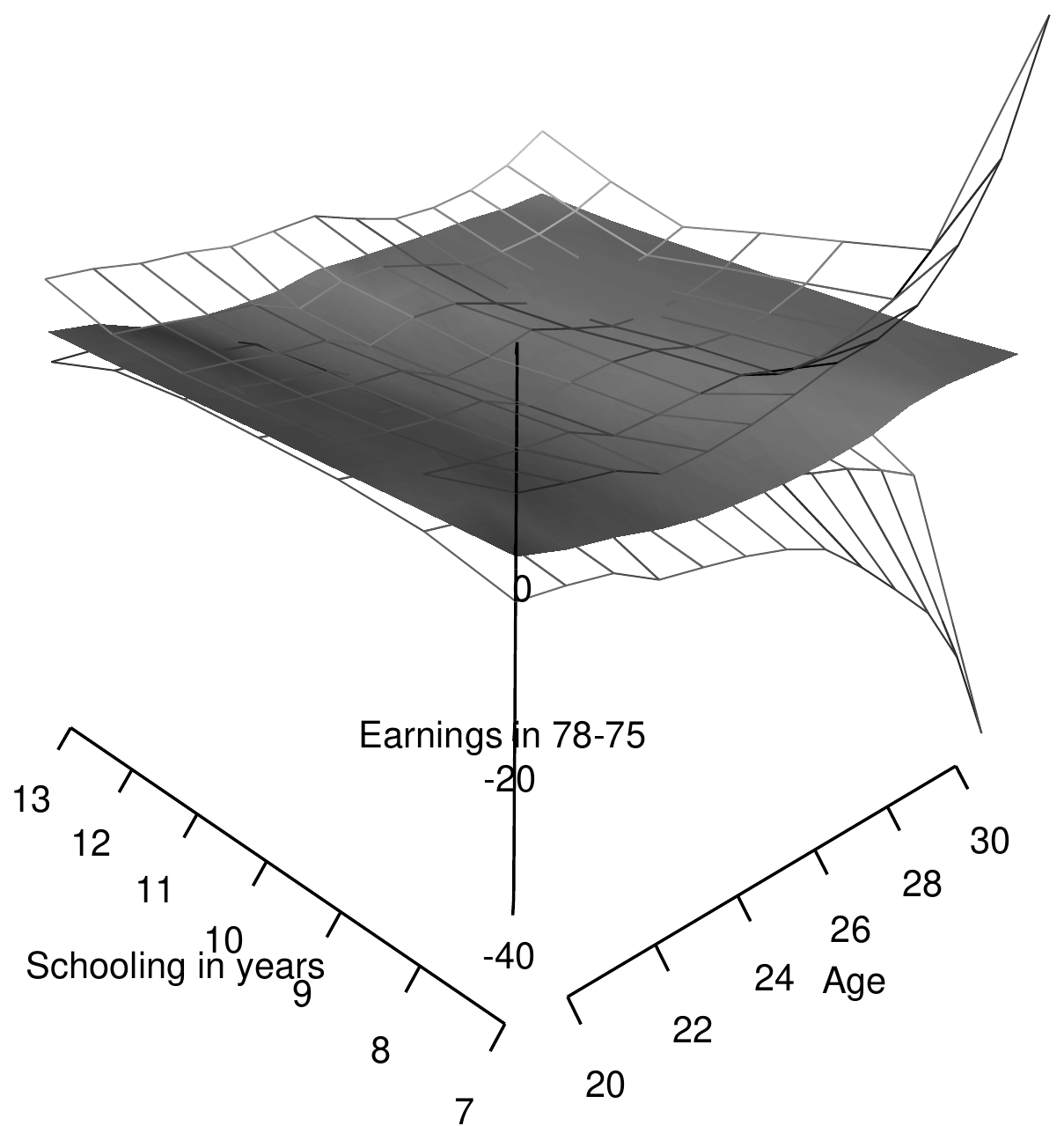}}
\subfigure[$\tau=20\%$]{\includegraphics[width=5cm, height = 5cm]{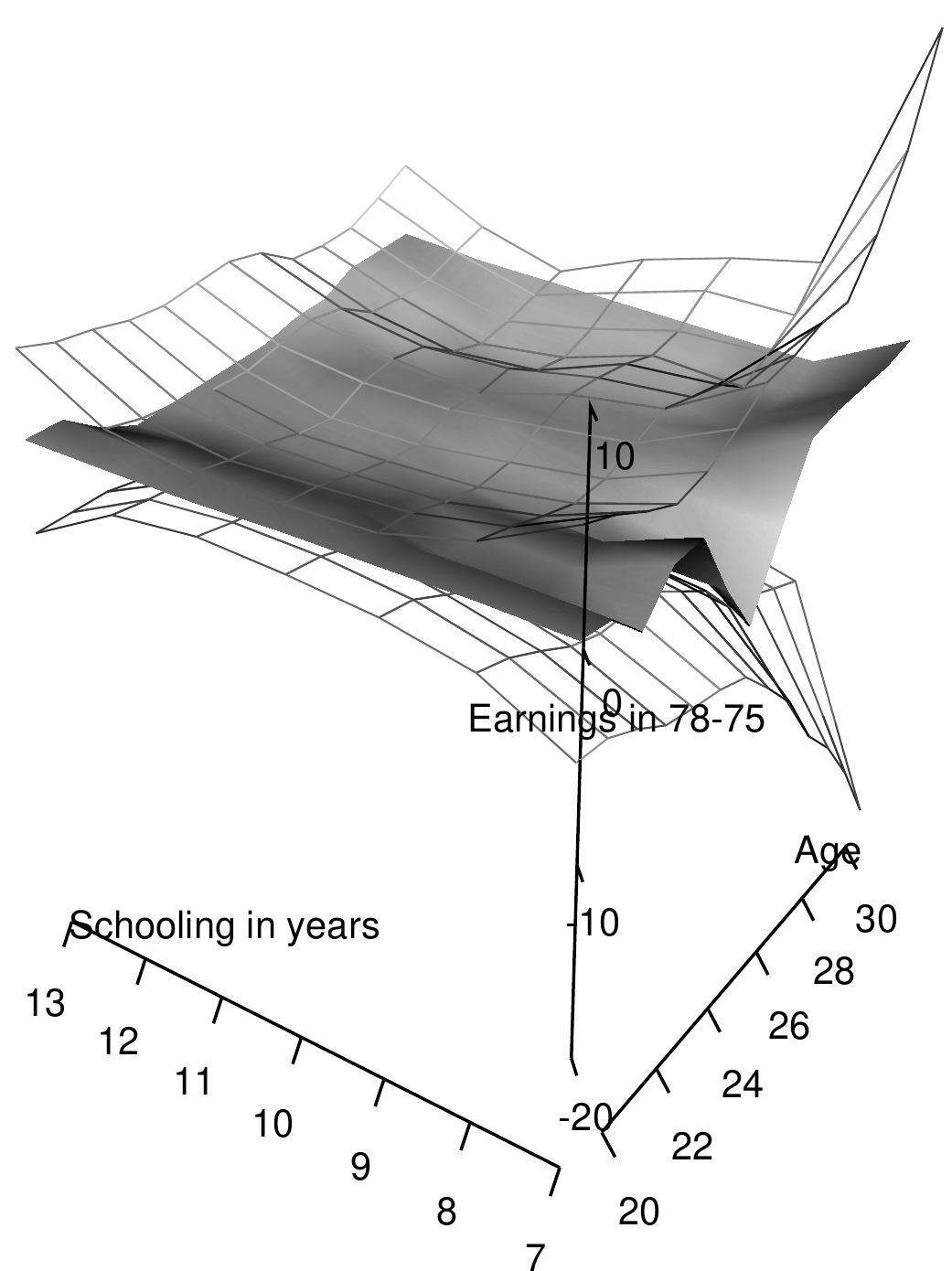}}
\subfigure[$\tau=30\%$]{\includegraphics[width=5cm, height = 5cm]{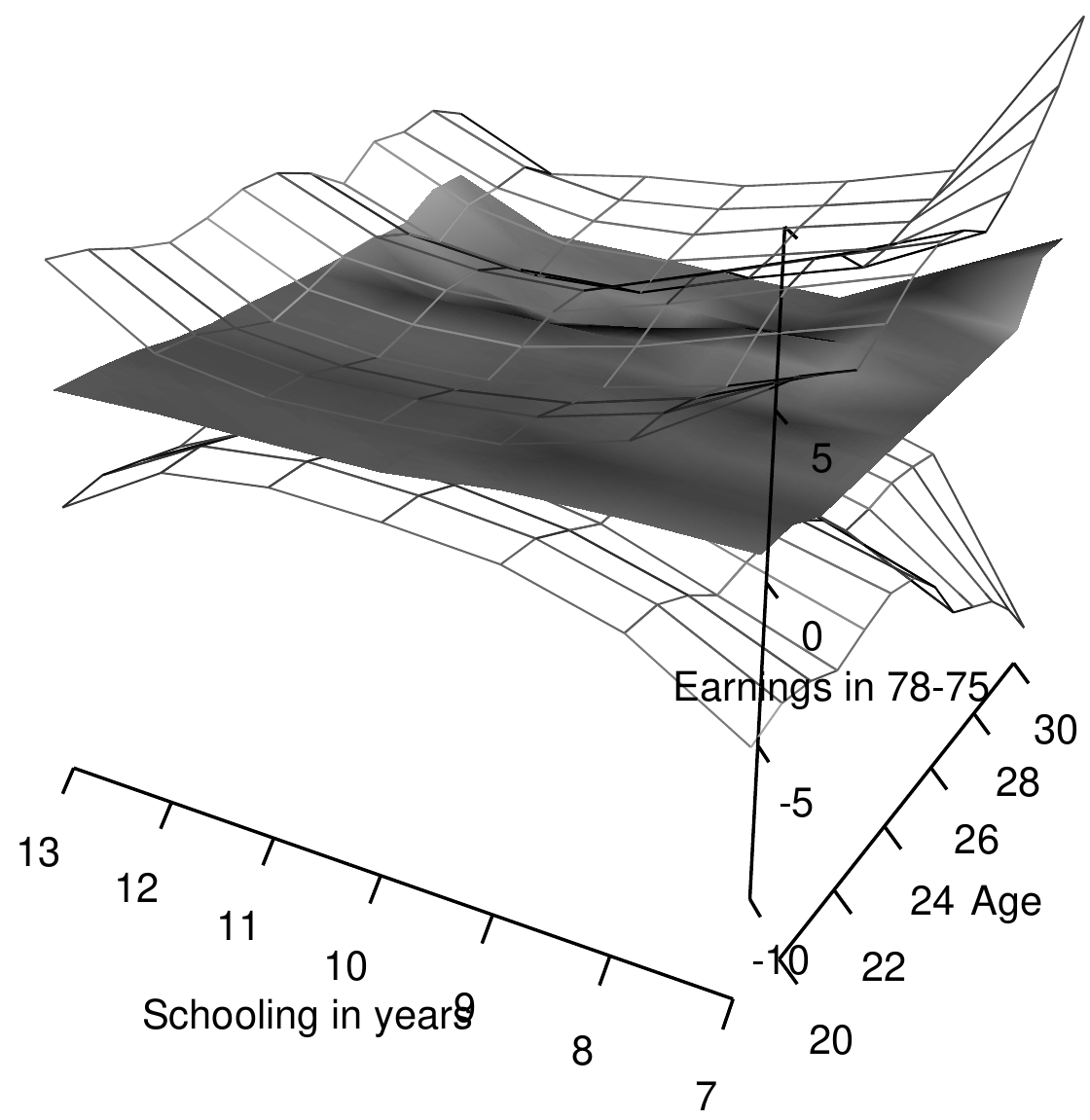}}
\subfigure[$\tau=50\%$]{\includegraphics[width=5cm, height = 5cm]{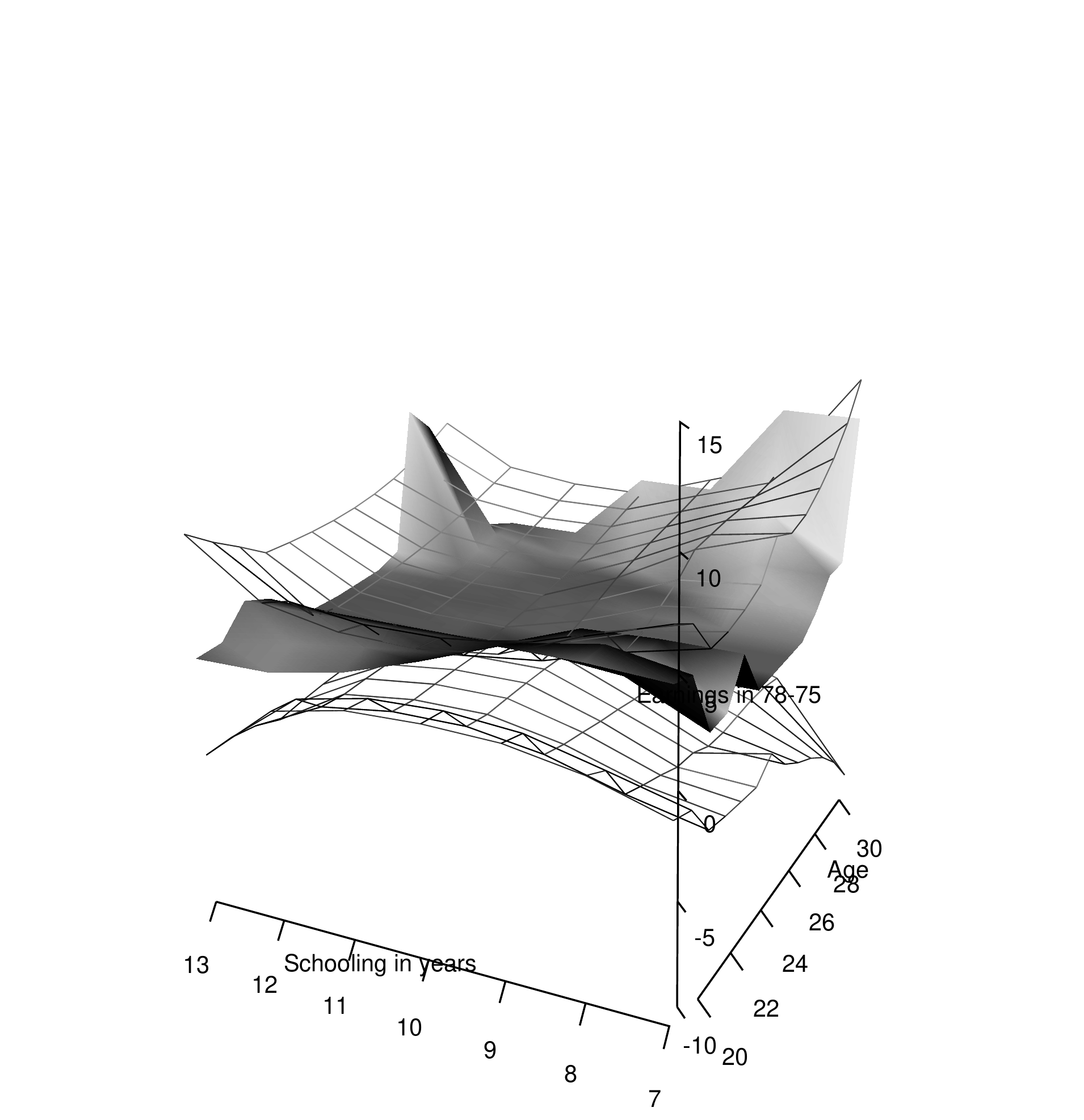}}
\subfigure[$\tau=70\%$]{\includegraphics[width=5cm, height = 5cm]{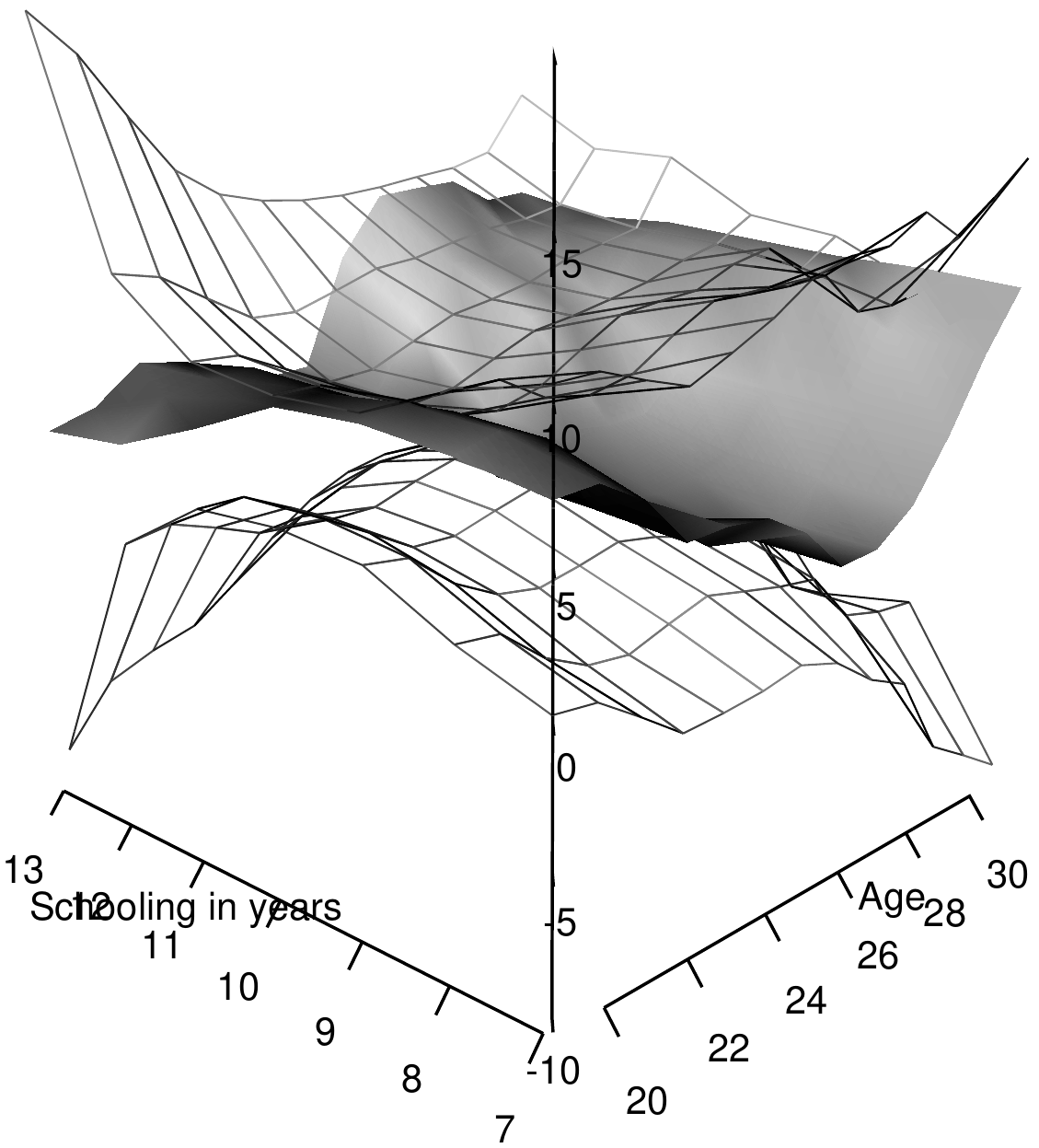}}
\subfigure[$\tau=80\%$]{\includegraphics[width=5cm, height = 5cm]{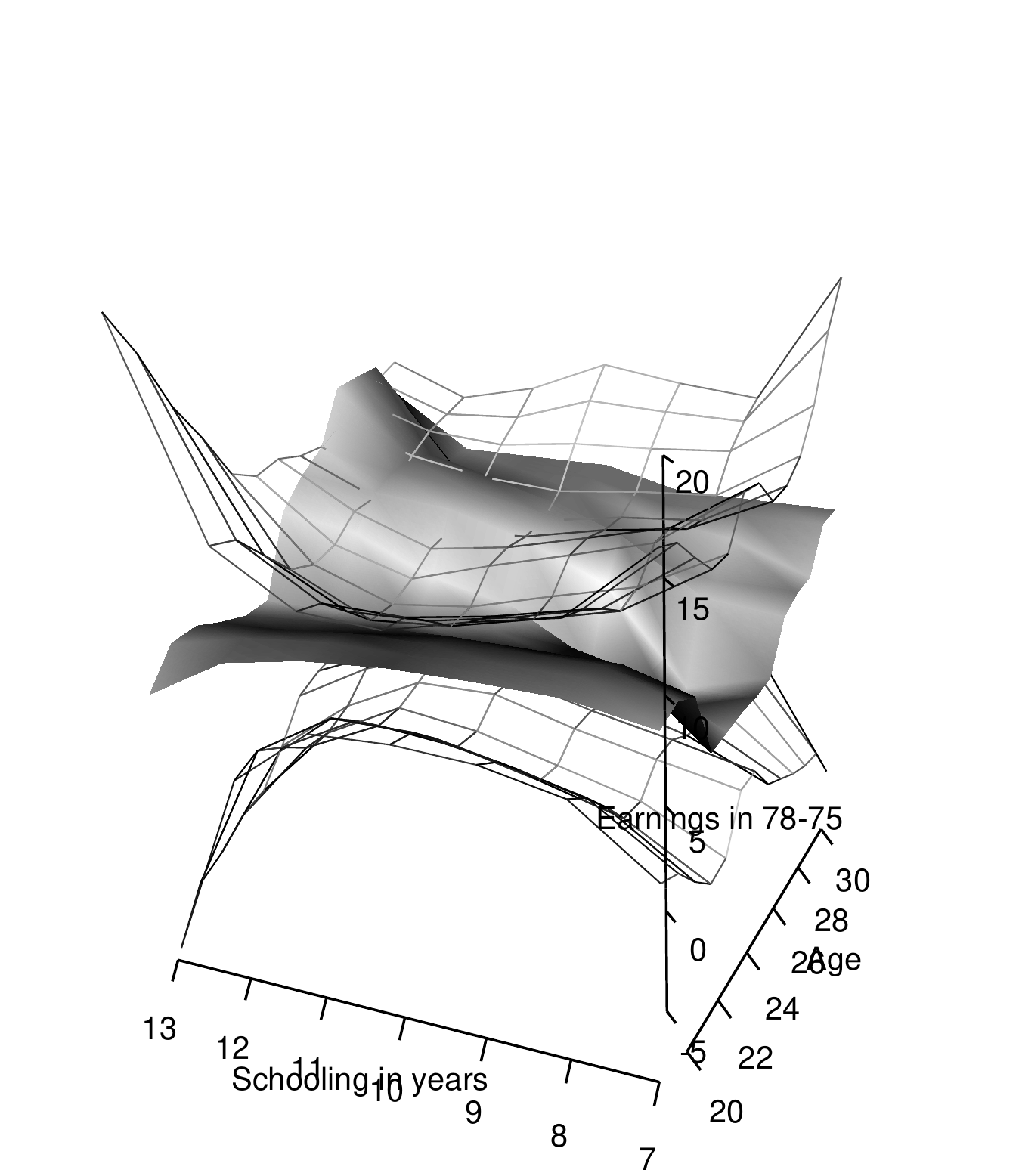}}
\subfigure[$\tau=90\%$]{\includegraphics[width=5cm, height = 5cm]{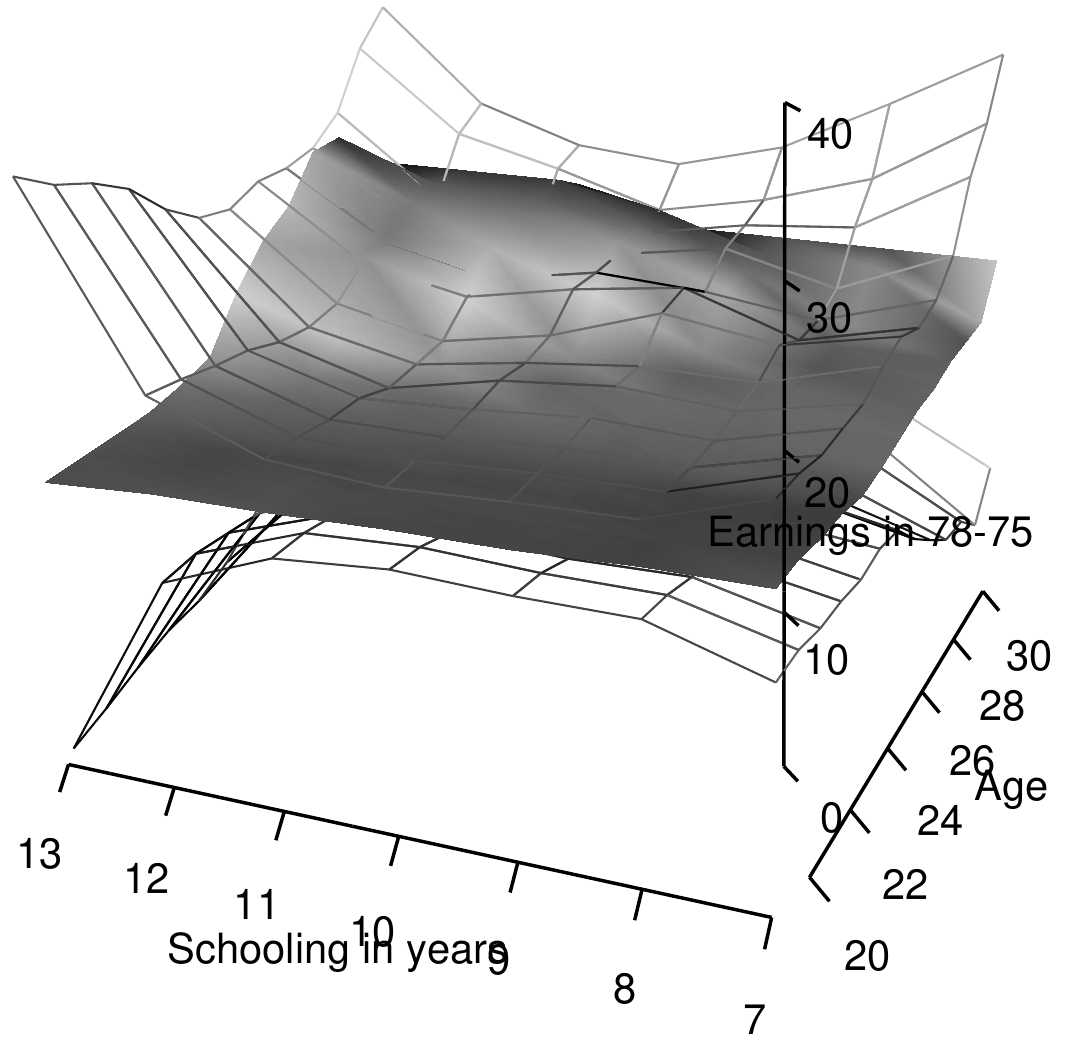}}
\caption{The conditional quantiles (solid surfaces) for the treatment group and the CCs (net surfaces) for the control group.}\label{lalonde.fvtr}
\end{figure}

The second observation obtained from comparing subfigures in Figure
\ref{lalonde.fvtr}, we find that the treatment has larger impact in
raising the upper bound of the earnings growth than improving the
lower bound. For lower quantile levels $\tau=10\%, 20\%$ and 30\%
the solid surfaces uniformly lie inside the CC of the control group,
while for $\tau=50\%, 70\%, 80\%$ and $90\%$, we see several
positive exceedances over the upper boundary of the CC of the
control group. Hence, the program tends to do better at raising the
upper bound of the earnings growth but does worse at improving the
lower bound of the earnings growth. In other words, the program
tends to increase the potential for high earnings growth but does
little in reducing the risk of negative earnings growth.

Our last conclusion comes from inspecting the shape of the surfaces:
conditioning on different levels of years of schooling (age), the
treatment effect is heterogeneous in age (years of schooling). The
most interesting cases occur when conditioning on high age and
high years of schooling. Indeed, when considering the cases of
$\tau=80\%$ and $90\%$, when conditioning on the years of schooling
at 12 (corresponding to finishing the high school), the earnings
increment of the treatment group rises above the upper boundary of
the CC of the control group. This suggests that the individuals who
are older and have more years of schooling tend to benefit more from
the treatment.

\section*{Supplementary Materials}
Section A contains the detailed proofs of Theorems 2.1, 2.3, 3.1 and
Lemmas 2.6 and 3.2, as well as intermediate results. Section B
contains some results obtained by other authors, which we use in our
study. We incorporate them here for the sake of completeness.



\bibliographystyle{dcu}
\bibliography{bibConfiBandQuant}
\appendix
\appendixpage
\section{Assumptions}\label{Sec:Ass}
\begin{enumerate}
\item[(A1)] $K$ is of order $s-1$ (see (A3)) has bounded support $[-A,A]^d$ for $A >0$ a positive real scalar, is continuously differentiable up to order $d$ with bounded derivatives, i.e.
$\partial^\sboldalpha K \in L^1(\R^d)$ exists and is continuous
for all multi-indices $\boldalpha \in \{0,1\}^d$
\item[(A2)] Let $a_n$ be an increasing sequence, $a_n\rightarrow\infty$ as $n\rightarrow\infty,$ and the marginal density $f_Y$ be such that
\begin{align}
(\log n) h^{-3d} \int_{|y|>a_n} f_Y(y)dy  = \CO(1) \label{A2}
\end{align}
and
$$
(\log n) h^{-d} \int_{|y|>a_n} f_{Y|\sbX}(y|\bx)dy  = \CO(1), \mbox{ for all }\bx \in \mathcal D
$$
as $n\rightarrow\infty$ hold.
\item[(A3)] The function $\theta_0(\bx)$ is continuously differentiable and is in H\"older class with order $s >d$.
\item[(A4)] $f_\sbX(\bx)$ is bounded, continuously differentiable and its gradient is uniformly bounded. Moreover, $\inf_{\sbx \in \mathcal D} f_\sbX(\bx) > 0$.
\item[(A5)] The joint probability density function $f(y,\bu)$ is bounded, positive and continuously differentiable up to $s$th order (needed for Rosenblatt transform).
The conditional density $f_{Y|\sbX}(y|\bx)$ exists and is
boudned and continuouly differentiable with respect to $\bx$.
Moreover, $\inf_{\sbx \in \mathcal D}
f_{Y|\sbX}\big(\theta_0(\bx)|\bx\big) >0$.
\item[(A6)] $h$ satisfies $\sqrt{nh^d} h^s \sqrt{\log n} \rightarrow 0$ (undersmoothing),
and $n h^{3d} (\log n)^{-2} \to \infty$.
\end{enumerate}
\begin{itemize}
\item[(EA2)]
$\sup_{\sbx\in\mathcal D}\left|\int v^{b_1} f_{\vep|\sbX}(v|\bx)
dv \right| <\infty$, for some $b_1>0$.
\end{itemize}
\begin{enumerate}
  \item[(B1)] $L$ is a Lipschitz, bounded, symmetric kernel. $G$ is Lipschitz continuous cdf, and $g$ is the derivative of $G$ and is also a density, which is Lipschitz continuous, bounded, symmetric and five times continuously differentiable kernel.
  \item[(B2)] $F_{\vep|\sbX}(v|\bx)$ is in $s'+1$ order H\"older class with respect to $v$ and continuous in $\bx$, $s' > \max\{2,d\}$. $\fx(\bx)$ is in second order H\"older class with respect to $\bx$ and $v$. $\E[\psi^2(\vep_i)|\bx]$ is second order continuously differentiable with respect to $\bx \in \mathcal D$.
\item[(B3)] $nh_0\brh^d \to \infty$, $h_0,\brh = \CO(n^{-\nu})$, where $\nu >0$.
\end{enumerate}

\begin{enumerate}
  \item[(C1)] There exist an increasing sequence $c_n$, $c_n\rightarrow\infty$ as $n\rightarrow\infty$ such that
\begin{align}
(\log n)^3 (nh^{6d})^{-1} \int_{|v|>c_n/2}  f_\vep(v)dv  = \CO(1) \label{C2},
\end{align}
as $n\rightarrow\infty$.
\item[(EC1)] $\sup_{\sbx\in\mathcal D}\left|\int v^b f_{\vep|\sbX}(v|\bx) dv
\right| <\infty$, for some $b > 0$.
\end{enumerate}
The assumptions (A1)-(A5) are assumptions frequently seen in the
papers of confidence corridors, such as \cite{WH:1989},
\cite{HS:2010} and \cite{GH:2012}. (EA2) and (EC1) essentially give
the uniform bound on the 2nd order tail variation, which is crucial
in the sequence of approximations for expectile regression.
(B1)-(B3) are similar to the assumptions listed in chapter 6.1 of
\cite{LR:2007}. (A6) characterizes the two conflicting conditions:
the undersmoothing of our estimator and the convergence of the
strong approximation. To make the condition hold, sometimes we need
large $s$ for high dimension, the smoothness of the true function.
(C1) and (EC1) are relevant to the theory of bootstrap, where we
need bounds on the tail probability and 2nd order variation.

\end{document}